\documentclass[fleqn,a4paper,10pt]{amsart}

\usepackage{graphicx}
\usepackage{tikz}
\usetikzlibrary{shapes,decorations.pathmorphing}
\usepackage{amsmath}
\usepackage{amssymb,amsthm}
\usepackage{comment}
\usepackage{hyperref}

\theoremstyle{definition}

\title[Eulerian triangulations:
two-point function and hull perimeter statistics]{Eulerian triangulations:
two-point function and hull perimeter statistics}
\author{Emmanuel Guitter}
\address{Institut de physique th\'eorique, Universit\'e Paris Saclay, CEA, CNRS, F-91191 Gif-sur-Yvette}
\email{emmanuel.guitter@cea.fr}

\begin{document}
\maketitle

\begin{abstract}
We present a new derivation of the distance-dependent two-point function for planar Eulerian triangulations
and give expressions for more refined generating functions where we also control hull perimeters. 
These results are obtained in the framework of a new recursion relation for slice generating functions
and extend similar results obtained recently for triangulations and quadrangulations.  
A number of explicit formulas are given for the statistics of hull perimeters in infinitely large random planar 
Eulerian triangulations. 
\end{abstract}

\section{Introduction}
\label{sec:introduction}

The enumeration of planar maps, namely connected graphs embedded on the sphere, has now become a classical field of investigation 
in combinatorics and probability theory. Among all families of maps, some are particularly simple and may serve as paradigms to test 
general ideas in the domain. 
This is the case of planar triangulations (i.e.\ maps whose all faces have degree $3$) or planar quadrangulations (all faces have degree $4$),
whose study gave rise to a lot of explicit results of all kinds, including many formulas for quite involved combinatorial quantities.

Another fundamental family of maps is that of \emph{planar Eulerian triangulations} which, so to say, consists of triangulations 
whose faces are colored alternatively in black and white. What makes this family particularly interesting is that maps with bi-colored faces 
provide in fact a more general framework in which many 
of the results for uncolored maps can be extended\footnote{Uncolored maps may in fact be viewed as a particular case of maps with bi-colored faces
by inflating their edges into faces of degree $2$ and coloring these new faces in black, the original faces being white.}
(see for instance \cite{BDG04}). Eulerian triangulations constitute the simplest family in this larger set and
their study also led to a number of explicit results. Eulerian triangulations are also dual of so-called bicubic maps, as first studied in \cite{TutteCPM},  and 
correspond to the simplest realization of so-called $3$-constellations \cite{BMS}.

Among the quantities amenable to an explicit formula is the so-called \emph{distance-dependent two-point function} which enumerates the maps at hand with a 
prescribed value of the \emph{graph distance} between two (randomly chosen) marked vertices. This two-point function can be given a particularly simple and elegant 
expression, both for triangulations and quadrangulations \cite{GEOD,BG12}, but also for Eulerian triangulations \cite{GEOD,AB11} provided, in this latter case, the distances
are measured along \emph{oriented paths} on the map (see below for a precise definition). 

Quite recently, a new approach was proposed to compute the distance-dependent two-point function, based on a direct decomposition of the underlying maps.
This approach was first implemented for triangulations in \cite{G15a}, then for quadrangulations in \cite{G15b}. It takes the form of a recursion relation for the generating function 
of so-called \emph{slices}, which are obtained by cutting and opening the maps at hand along some particular geodesic path. One advantage of the method is that it
not only controls the distance between the two marked vertices on the map, but also controls 
 the associated \emph{hull perimeter at distance $d$} which, so to say, is the length of the closed curve separating the two marked vertices and lying at a fixed distance 
$d$ from the first vertex (see \cite{Krikun03,Krikun05,CLG14a,CLG14b,G16}). 
A number of new results on the statistics of hull perimeters, depending on both the distance between the two marked vertices and the
 distance $d$ at which the perimeters are measured, were therefore obtained in \cite{G16}, where a particular emphasis was put on the limit of 
 infinitely large maps.   
 
 The purpose of this paper is to extend the approach of \cite{G15a,G15b} to the case of planar Eulerian triangulations, so as to again provide a new, more direct and constructive 
 derivation of their distance-dependent two-point function and, in a second step, to compute some statistical properties of their hull perimeters.
 
The paper is organized as follows: Section \ref{sec:eulerian} recalls a number of basic definitions about Eulerian triangulations
(Section \ref{sec:basics}) and discusses the connection between the distance-dependent two-point function and the generating function of slices (Section \ref{sec:slices}).
We then remind the reader of a classical relation satisfied by these slice generating functions (Section \ref{sec:recurold}), whose derivation is recalled in Appendix A.
Section \ref{sec:recurnew} explains how slices may be decomposed in a recursive way similar to that
of \cite{G15b} (Section \ref{sec:decomp}) and discusses the resulting recursion relation for slice generating functions as well as the particular form of its \emph{kernel} 
$\mathcal{K}$ (Section \ref{sec:kernel}). This kernel itself involves a particular generating function $\Phi$, which enumerates maps with a boundary 
and with a number of forbidden boundary-boundary connections (Section \ref{sec:h2i}).
The explicit computation of $\Phi$ is presented in Section \ref{sec:phi}, where we first explain how to write a closed system involving
$\Phi$ and another similar generating function $\Omega$ (Section \ref{sec:system}) and then show how to solve this system (Section \ref{sec:solsystem}).
The explicit knowledge of $\Phi$ allows us to rewrite our recursion relation in a form which can then be solved in a straightforward way (Section \ref{sec:solrecur}).
This eventually allows us to get the desired explicit expression for the distance-dependent two-point function of planar Eulerian triangulations.
Section \ref{sec:hull} is devoted to a study of the hull perimeter statistics. After some basic definition (Section \ref{sec:perimeter}), we derive an explicit formula 
for a particular generating function which generalizes the two-point function and offers some additional control on hull perimeters (Section \ref{sec:genhp}). We finally use this expression to
obtain a number of explicit results on the statistics of hull perimeters in the limit of large planar Eulerian triangulations (Section \ref{sec:stathp}).
We gather a few concluding remarks in Section \ref{sec:conclusion}.

\section{Eulerian triangulations, two-point function and slices}
\label{sec:eulerian}
\subsection{Eulerian triangulations: basic facts}
\label{sec:basics}
A \emph{planar Eulerian triangulation} is a planar map whose all faces have degree $3$ and are colored in black
or white in such a way that any edge of the map is incident to a face of each color. Edges in a planar Eulerian triangulation are naturally oriented by 
demanding that the face on their left be black. Otherwise stated, edges are oriented counterclockwise around black faces and 
clockwise around white faces.

It is easily seen that the length of any oriented cycle in the map is a multiple of $3$ so that we may assign colors $0$, $1$ or $2$
to the vertices in such a way that the color increases by $1\!\mod 3$ when we follow an oriented edge. In other words, 
the colors appear cyclically in the order $0\to 1\to 2$ counterclockwise around black faces and clockwise around 
white faces\footnote{The assignment 
of colors is moreover unique, up to a global cyclic permutation.}.
It is also easily seen that any two vertices are accessible from one another by following some oriented path. 

A planar Eulerian triangulation is \emph{pointed} if it has a marked vertex $v_0$, called the \emph{origin}. To each vertex $v$ of the map, we may then assign 
its ``oriented distance" $d(v)$ (later called ``distance" for short) defined as the length of any shortest oriented path from $v_0$ to $v$. An oriented edge of the map
is said of type $(\ell,\ell')$ if it points from a vertex $v$ with $d(v)=\ell$ to a vertex $v'$ with $d(v')=\ell'$.
Clearly, if we give the color $0$ to the origin $v_0$, the distance $d(v)$ of any vertex $v$ satisfies $d(v)=c(v)\!\mod 3$. 
In particular, the distance of vertices increases by $1\!\mod 3$ when we follow an oriented edge. Since distances 
around a face vary by at most $2$,  oriented edge come in two flavors only:
\emph{short-edges} of type $(\ell,\ell+1)$ or \emph{long-edges} of type $(\ell+2,\ell)$, for some $\ell\geq 0$.

\begin{figure}
\begin{center}
\includegraphics[width=5cm]{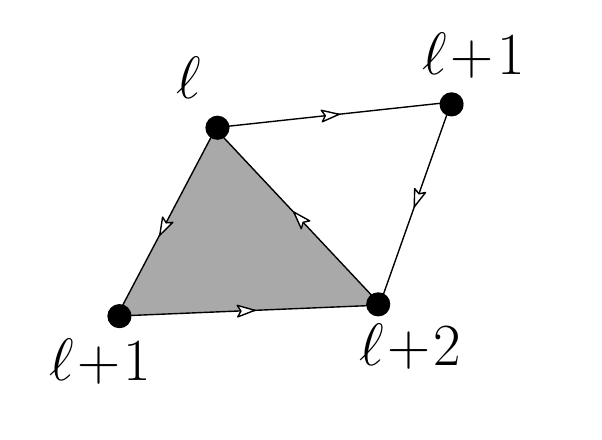}
\end{center}
\caption{The configuration of distances at the vertices around two adjacent white and black (represented here in gray) faces sharing a long-edge.
The two faces are of the same type $\ell$.}
\label{fig:blackwhite}
\end{figure}
Any face in the map is incident to exactly one long-edge and two short-edges of type  $(\ell+2,\ell)$, $(\ell,\ell+1)$ and $(\ell+1,\ell+2)$ respectively, 
for some $\ell \geq 0$: the three edges appear in this order counterclockwise around the face if it is black and clockwise if it is white. 
Such faces will be refereed to as black and white faces of type $\ell$.
In particular, any long-edge $(\ell+2,\ell)$ separates a black and a white face of \emph{the same type $\ell$} \footnote{This is not necessarily true for a short-edge 
which may separate black and white faces whose types differ by $1$.} (see figure \ref{fig:blackwhite}). As a consequence, a pointed planar Eulerian triangulation has the same number
of black faces and white faces, which is also its number of long-edges. Note finally that \emph{erasing all its long-edges} provides a canonical way to transform   
a pointed planar Eulerian triangulation into a particular pointed planar quadrangulation, i.e.\ a map whose all faces have degree $4$ \footnote{Not all 
pointed planar quadrangulations however are images of pointed Eulerian triangulations by this transformation.}.

\subsection{Two-point function and slice generating function}
\label{sec:slices}
\begin{figure}
\begin{center}
\includegraphics[width=15cm]{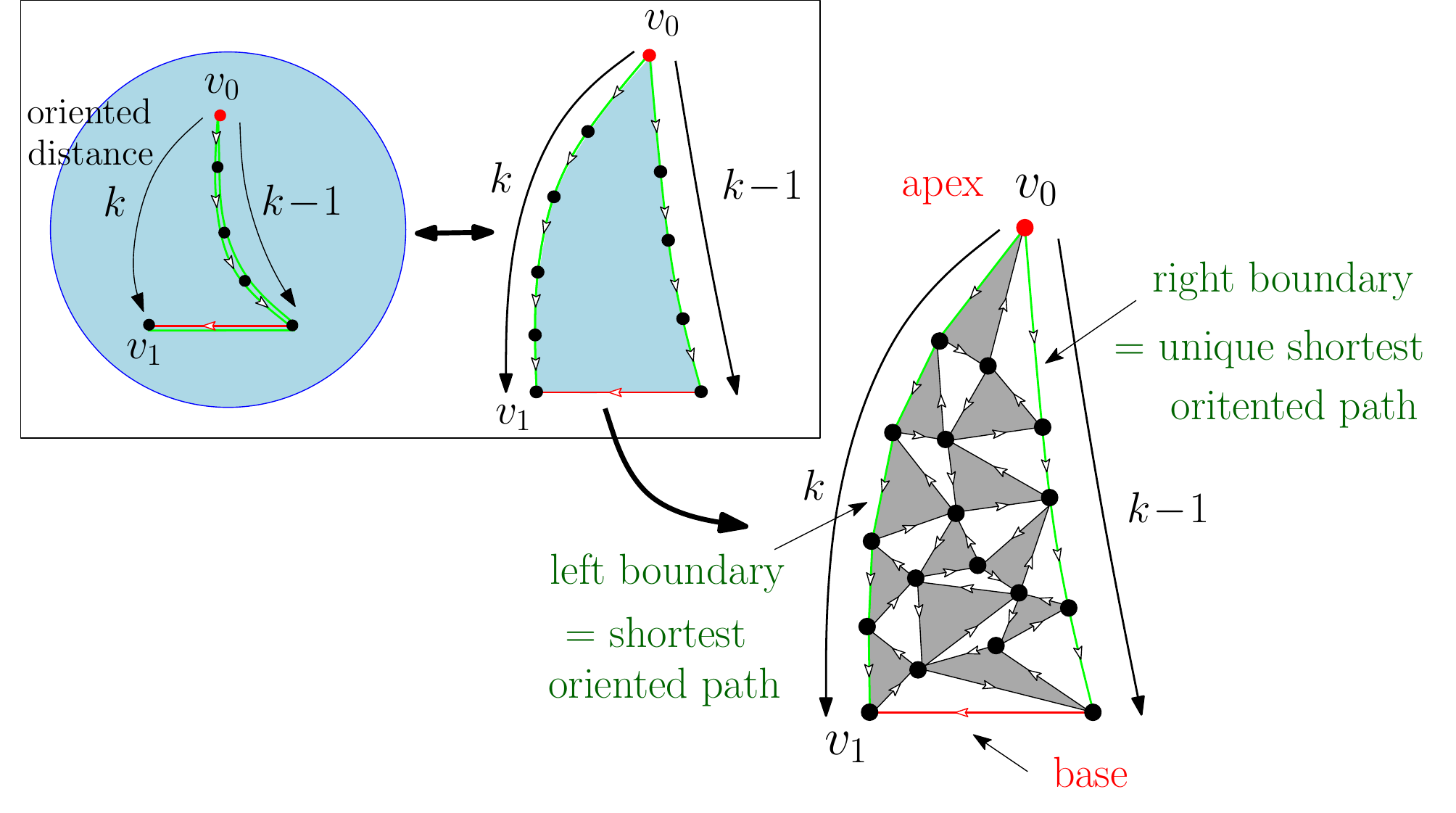}
\end{center}
\caption{In the frame: a schematic picture of the transformation from a pointed planar Eulerian triangulation (here represented on the $2$-dimensional sphere,
the light-blue background being supposedly filled with black and white faces) 
with a marked oriented edge of type $(k-1,k)$ into a $k$-slice by cutting along the leftmost
backward shortest path (in green) from the endpoint $v_1$ of the marked oriented edge (red segment) to the origin $v_0$ (red dot).
The obtained slit is then opened so as to form the outer face of the slice.
The origin $v_0$ becomes the apex of the $k$-slice and the marked oriented edge its base. On the right: a more detailed example of 
a $k$-slice where we now represented the black and white faces and the orientation of edges.}
\label{fig:slice}
\end{figure}
In this paper, we adopt a simple definition of the \emph{distance-dependent two-point function} $G_k\equiv G_k(g)$ as the generating function of 
pointed planar Eulerian triangulations \emph{with a marked oriented edge of type $(k-1,k)$} for some $k\geq 1$, enumerated with 
a weight $g$ per white face (recall that the number of white faces is also that of black faces as well as that of long-edges).
In other words, maps enumerated by $G_k$ are endowed with two marked vertices at oriented distance $k$ from one another, the origin $v_0$ of
the map and the endpoint $v_1$ of the marked edge, but the data of these two vertices is supplemented by the precise choice of the marked edge itself, i.e.\ by
the data of some particular marked oriented edge pointing towards $v_1$ and whose first extremity lies at distance $k-1$ from $v_0$ (note that
such an oriented edge always exists for any $v_1$ at distance $k$ from $v_0$). This definition is intended to eventually lead to a simple expression for the two-point function $G_k$.

The maps enumerated by $G_k$ may be transformed into so-called \emph{$k$-slices} by a simple (reversible) cutting procedure (see figure \ref{fig:slice} for an illustration). 
The transformation is as follows: starting from the endpoint $v_1$ of the marked edge (at distance $k$
from the origin $v_0$), we first follow the marked edge backwards to its other extremity (at distance $k-1$), then follow the
\emph{leftmost} edge leading backwards to a vertex at distance $k-2$ and continue the process by following backwards leftmost
edges to vertices at distances $k-3, k-4, \cdots$ until we eventually reach the origin $v_0$. The followed path constitutes what we may call the \emph{leftmost
backward shortest path} from $v_1$ to $v_0$. Cutting the map along this path and opening it results into
is a particular planar map, which we call a $k$-slice, fully characterized by the following properties (see figure \ref{fig:slice}):
\begin{itemize}
\item[(s1)] A $k$-slice is a planar \emph{rooted} map, i.e.\ a map with a marked oriented edge, later called the \emph{base}, and whose \emph{outer face}, i.e.\ 
the face lying on the left of the base, has degree $2k$. The boundary of this outer face is moreover a simple closed curve (which corresponds, after cutting of the original Eulerian map enumerated by $G_k$, to two copies of the leftmost 
backward shortest path) made of three parts parts:
(i) the base itself (which corresponds to the marked oriented edge of type $(k-1,k)$ in the original Eulerian map), 
(ii) a left boundary made of $k$ oriented edges leading from an \emph{apex} $v_0$ (which corresponds to the origin in the original Eulerian map) to the endpoint $v_1$ of
the base and (iii) a right boundary made of $k-1$ oriented edges leading from $v_0$ to the first extremity of the base\footnote{The ``left" and ``right" denominations
for the boundaries correspond to the usual terminology which, regardless of orientations, refers to the position of the boundaries with respect to the bulk of the slice in a picture where the latter 
is represented with its base at the bottom, as in figure \ref{fig:slice}. As an oriented path however, the left (respectively right) boundary is on the right 
(respectively on the left) of the bulk of the slice, i.e.\ it has the bulk on its left (respectively on its right).}.

\item[(s2)] The inner faces, i.e.\ the faces other than that the outer face, all have degree $3$ and form what we call the bulk of the $k$-slice.
These faces are colored in black and white so that any inner edge (i.e.\ an edge incident to inner faces only) is incident to one face of each color.
As before, these edges are oriented with their black incident face on the left.  As for the (already oriented) edges incident to the outer face, those of the left boundary 
are incident to black inner faces only on their left while the root edge and the edges of the right boundary are incident to white inner faces only 
on their right.

\item[(s3)] The left boundary of the $k$-slice is a shortest oriented path from $v_0$ to $v_1$ within the $k$-slice. The right boundary of the $k$-slice is the \emph{unique} 
shortest path from $v_0$ to the origin of the root edge within the $k$-slice.
\end{itemize}
This last property is a direct consequence of our choice of cutting along the leftmost backward shortest path. The above properties (s1)--(s3) define what we 
call $k$-slices.
For convenience, the "single-edge-map", i.e.\ the map reduced to a single edge oriented from $v_0$ to $v_1$ and an outer face of degree $2$ is considered a $1$-slice.

Clearly, our cutting procedure transforms any pointed planar Eulerian triangulation enumerated by $G_k$ into a $k$-slice. Conversely, any $k$-slice but the single-edge-map encodes a particular map enumerated by $G_k$, the latter being easily obtained by re-gluing the left boundary of the slice (of length $k$) to the right boundary supplemented by the 
base (of total length $k-1+1=k$). Clearly, all the faces, vertices and edges of the original Eulerian triangulation are recovered in the $k$-slices (the vertices and edges along the
leftmost backward shortest path appearing twice in the slice) and the oriented distance $d(v)$ from the origin $v_0$ to some vertex $v$ in
the original Eulerian triangulation is simply identified as the oriented distance from the apex $v_0$ to the corresponding vertex in the $k$-slice. 
In particular, each inner face in the $k$-slice is, as before, a black or a white face of type $\ell$, i.e.\ incident to exactly one long-edge and two short-edges 
of type  $(\ell+2,\ell)$, $(\ell,\ell+1)$ and $(\ell+1,\ell+2)$ respectively, for some $\ell \geq 0$.

We call $R_k= R_k(g)$ ($k \geq 1$) the generating function of $\ell$-slices \emph{with $1 \leq \ell \leq k$}, enumerated with a weight $g$ per inner face.
Note that the "single-edge-map", which has no inner face, contributes a term $1$ to $R_k$ for all $k\geq 1$. From the above transformation, we have 
the correspondence $R_k=1+\sum_{\ell=1}^k G_\ell$ (with a first term $1$ corresponding to the contribution of the single-edge-map) or, by inversion
\begin{equation*}
G_k=R_k-R_{k-1}-\delta_{k,1}\ ,
\end{equation*} 
with the convention $R_0=0$. Computing the distance-dependent two-point function $G_k$ of planar Eulerian triangulation therefore reduces to computing 
the corresponding $k$-slice generating function $R_k$.

\subsection{A classical relation for slice generating functions}
\label{sec:recurold}
The question of computing $R_k$ in the context of planar Eulerian triangulations was first addressed in \cite{GEOD}, where it was shown that
 $R_k$ satisfies the following relation:
\begin{equation}
R_k=1+g\, R_k(R_{k+1}+R_{k-1})
\label{eq:recurold}
\end{equation}
for $k\geq 1$, with $R_0=0$. This relation follows from a simple decomposition of $k$-slices into 
appropriate sub-slices and its derivation is recalled in Appendix A below. 
Note that it fixes all the $R_k$'s order by order in $g$ by demanding that $R_k=1+O(g)$ for all $k\geq1$, as 
required by the definition of $R_k$. At this stage, it is important to notice that the
explicit form of $R_k$ given in  \cite{GEOD} was not derived strictly speaking from the relation itself but 
was simply the result of some educated guess (the solution being unique, it is enough to check that the 
proposed guess actually solves the equation and satisfies the small $g$ expansion requirement). In this
respect, the approach of  \cite{GEOD} is not fully satisfactory.
Later, the expression of $R_k$ was recovered, now in a constructive way, in  \cite{AB11}, where it was shown that the
$R_k$'s are coefficients in a particular \emph{multi-continued fraction} expansion for some more
global generating function.  

In the next Section of this paper, we will present a new constructive approach to compute $R_k$,
based on a more direct recursive relation linking $R_k$ to $R_{k-1}$. By a slight refinement, the method 
will then allow us to explore properties 
of hull perimeters in planar Eulerian triangulations.

When $k\to\infty$, $R_k$ tends to a limit $R_\infty$ which is the generating function of $\ell$-slices with arbitrary $\ell\geq 1$.
From \eqref{eq:recurold}, $R_\infty$ is determined by the equation
\begin{equation*}
R_\infty=1+2 g\, R_\infty^2\ ,
\end{equation*}
with the requirement that $R_\infty=1+O(g)$, hence 
\begin{equation*}
R_\infty=\frac{1-\sqrt{1-8g}}{4g}\ .
\end{equation*}

\begin{figure}
\begin{center}
\includegraphics[width=4cm]{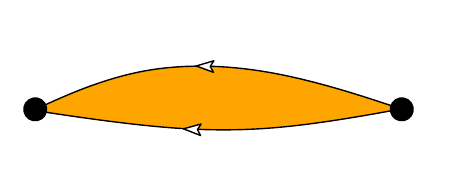}
\end{center}
\caption{A schematic picture of a bundle (the orange background is supposedly filled with black and white faces).
The generating function for bundles is $R_1$.}
\label{fig:bundle}
\end{figure}
Finally, for $k=1$, $R_1$ enumerates $1$-slices, i.e.\ slices for which the apex is identical to the first
extremity of the base. This corresponds to maps with a boundary of length $2$ which, when the map is not
reduced to the single-edge-map, are formed of two distinct edges both pointing from the first extremity
of the base to its endpoint (see figure \ref{fig:bundle} for an illustration) and a number of black or white inner faces in-between. 
Such maps with a boundary of length $2$ we be called \emph{bundles} in the following.

\section{A new direct recursion relation for slice generating functions}
\label{sec:recurnew}
We now come to the main point of this paper, namely the derivation of a new recursion relation for $R_k$,
or, more precisely, for 
\begin{equation*}
T_k\equiv R_k-R_1\ ,
\end{equation*}
$k\geq 1$, which is the generating function for $\ell$-slices with $2\leq \ell \leq k$ (note that, in particular, $T_1=0$). 
At this stage, we simply adapt to our problem of Eulerian triangulations a construction developed in \cite{G15b} for quadrangulations.
As in \cite{G15b}, the recursion relation takes the form
 \begin{equation}
T_k=\frac{R_1^2 (T_{k-1}+R_1)\, \Phi(T_{k-1})}{1-R_1(T_{k-1}+R_1)\, \Phi(T_{k-1})}\ , \qquad \Phi(T)\equiv \Phi(T,g)=\sum_{i\geq 2} h_{2i}(g)\, T^{i-2}\ ,
\label{eq:newrecur}
\end{equation}
where $h_{2i}(g)$, $i\geq 2$ enumerate particular maps whose definition is detailed below.
This relation is the result of some splitting of the $\ell$-slices into various domains upon cutting along particular paths, as we describe now.

\subsection{The slice decomposition}
\label{sec:decomp}
\begin{figure}
\begin{center}
\includegraphics[width=8cm]{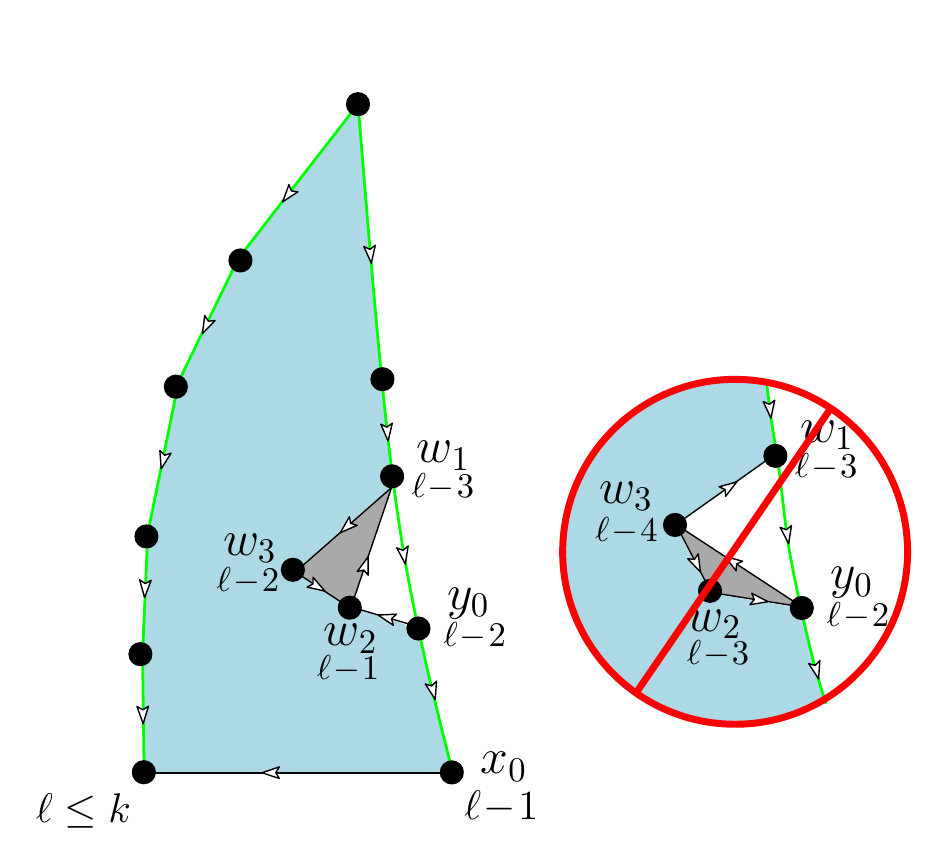}
\end{center}
\caption{Configuration of distances in the vicinity of the right-boundary vertex $y_0$ at distance $\ell-2$ from the apex in an $\ell$-slice enumerated by $T_k$
(see text -- the light-blue background is supposedly filled with black and white faces). The configuration displayed on the right is forbidden.}
\label{fig:dividingzero}
\end{figure}
Our decomposition is similar to that of \cite{G15b}. 
We start with an $\ell$-slice with $2\leq \ell \leq k$, as enumerated by $T_k$, and select on the slice some particular \emph{dividing line} defined as follows:
assuming $\ell\geq 3$, we look at the white face on the right of the (unique) oriented right-boundary edge of type $(\ell-3,\ell-2)$.
A priori, the long-edge incident to this face either starts from the right-boundary vertex $y_0$ at distance $\ell-2$ or ends
at the right-boundary vertex $w_1$ at distance $\ell-3$ (see figure \ref{fig:dividingzero} for an illustration). However, it is easily seen that  
this long-edge cannot start from $y_0$ as otherwise, the white face
would be of type $\ell-4$ (requiring in particular $\ell\geq 4$), and so would be the black face on the left of the long-edge. This would imply that $y_0$ is adjacent, via an edge in the bulk, 
to a vertex $w_2$ at distance $\ell-3$ from the apex (see figure \ref{fig:dividingzero}), in contradiction with (s3) which states that the right boundary of the $\ell$-slice is the unique
shortest path between the apex and the first extremity $x_0$ of the base. The long-edge therefore ends at $w_1$ and starts at a vertex $w_2$ 
at distance $\ell-1$ from the apex. The white face is therefore of type $\ell-3$ and so is the black face on the other side
of the long-edge. This black face, incident to $w_1$ and $w_2$, is therefore also incident to a third vertex $w_3$ at distance 
$\ell-2$. Note that $w_3$ cannot be equal to $y_0$ as otherwise, the edge from $w_1$ to $w_3$ would connect $w_1$ to $y_0$ within the bulk 
and create a second shortest path between the apex and the first extremity $x_0$ of the base, in contradiction with  (s3). As for $w_2$, since it is adjacent to $w_3$, 
it cannot be equal to the first extremity $x_0$ of the base since, because of (s3), $x_0$ cannot be adjacent to a vertex at distance $\ell-2$ 
and different from $y_0$. We thus have a \emph{two-step path} (i.e.\ a sequence of two consecutive edges) 
$y_0\to w_2\to w_3$ starting from $y_0$ satisfying $d(w_2)=\ell-1$ and $d(w_3)=\ell-2$ with
$w_2\neq x_0$ and $w_3\neq y_0$. We may now pick the \emph{leftmost} such two-step path $y_0\to x_1 \to y_1$ starting from $y_0$, 
where $d(x_1)=\ell-1$ and $d(y_1)=\ell-2$ with $x_1\neq x_0$ and $y_1\neq y_0$.

\begin{figure}
\begin{center}
\includegraphics[width=11cm]{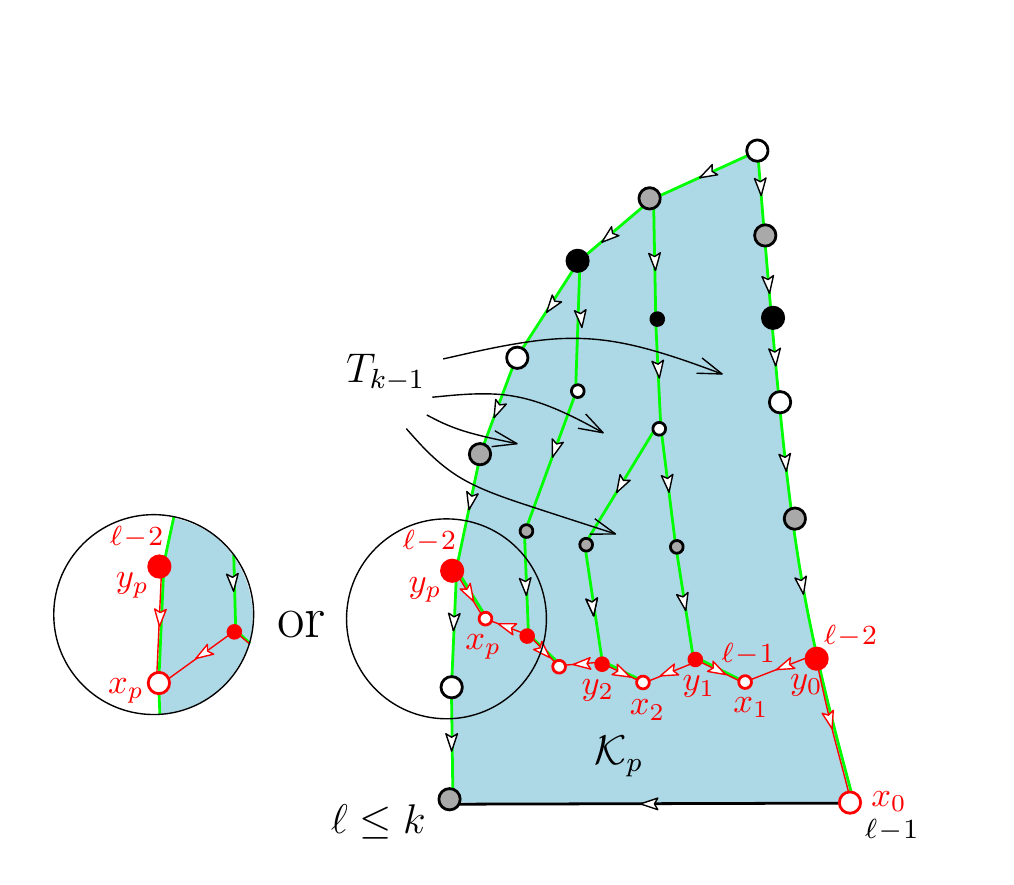}
\end{center}
\caption{A schematic picture of the construction of the dividing line (in red) in an $\ell$-slice enumerated by $T_k$ (the light-blue background is supposedly filled with black and white faces).
Vertices along the dividing line, of total length $2p+1$, alternate between vertices $x_m$ at distance $\ell-1$ (open circles) and vertices $y_{m'}$ at $\ell-2$ (filled circles). 
The dividing line hits the left boundary either at $x_p$ (distance $\ell-1$) or at $y_p$ (distance $\ell-2$). Drawing leftmost backward shortest paths from the vertices $x_m$ to the apex,
we see that to each two-step path $y_{m-1}\to x_{m}\to y_{m}$ of the dividing line is associated
a sub-slice (with base $y_{m-1}\to x_m$) in the upper part. Each of these sub-slices is enumerated by $T_{k-1}$ while the lower part is enumerated by $\mathcal{K}_p$.}
\label{fig:dividingone}
\end{figure}

This construction may be repeated as follows: we first draw the leftmost
backward shortest path from $x_1$ to the apex, taking the (backwards oriented) edge $x_1\to y_1$ as first step. This path serves as a right boundary
\emph{for the part of the slice lying in-between the path itself and the left boundary of the slice} and is the unique shortest path from the apex to $x_1$ in this domain. We again pick 
the leftmost two-step path $y_1\to x_2 \to y_2$ starting from $y_1$, 
where $d(x_2)=\ell-1$ and $d(y_2)=\ell-2$ with $x_2\neq x_1$ and $y_2\neq y_1$ (that such a path exists is proved by exactly the same arguments as above).
Continuing the process, it is easily shown that the obtained path $x_0\to y_0 \to x_1\to y_1\to x_2 \to y_2 \to x_3 \to y_3 \cdots$ cannot form loops
(see \cite{G15b} for a detailed argument) and is a simple path which eventually ends by reaching the left boundary after $p$ steps, either at the (unique) left-boundary vertex 
$x_p$ at distance $\ell-1$ 
from the apex or at the (unique) left-boundary vertex $y_p$ at distance $\ell-2$ (see figure \ref{fig:dividingone} for an illustration). 
As in \cite{G15b}, we will call this path the dividing line, as it separates the slice into two domains: an upper part containing the apex and a lower part
containing the base. 
\begin{figure}
\begin{center}
\includegraphics[width=9cm]{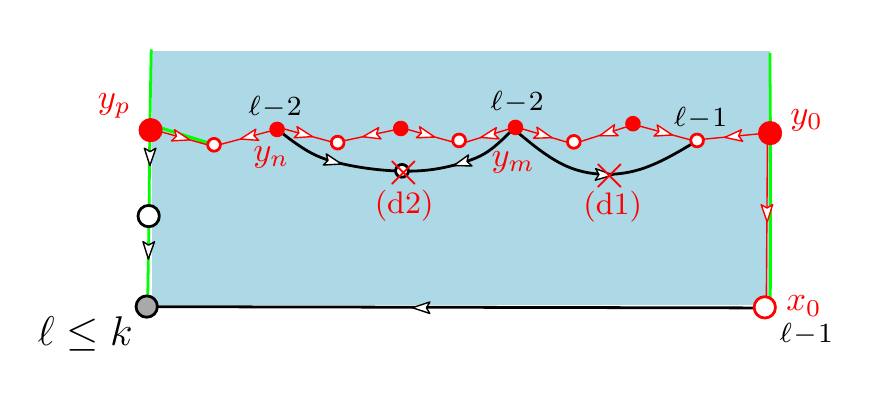}
\end{center}
\caption{An illustration of the two constraints of properties (d1) and (d2) fulfilled by the dividing line (see text).}
\label{fig:forbidden}
\end{figure}
By construction the dividing line, of length $2p+1$, satisfies the following two properties (see figure \ref{fig:forbidden} for an illustration):
\begin{itemize}
\item[(d1)] Two vertices of the dividing line cannot be linked by an edge lying strictly in the lower part\footnote{Note that this property is non-trivial only
if the two vertices are at respective distances $\ell-1$ and $\ell-2$ as otherwise, having a direct edge between them is already trivially forbidden by congruence.}.
\item[(d2)] Two distinct vertices $y_m$ and $y_n$ ($n\neq m$) of the dividing line, lying at distance $\ell-2$ from the apex, cannot have a common neighbor strictly in the lower part.
\end{itemize}
Both properties follow from the fact that, in the construction of the dividing line, we always picked the leftmost two-step paths
$y_m\to x_{m+1} \to y_{m+1}$. Note that two vertices $x_m$ and $x_n$ of the dividing line, lying now at distance $\ell-1$ from the apex,
may have a common neighbor strictly below the line (this common neighbor is then at distance $\ell$ and this does not lead to
any contradiction).

If $\ell=2$, $y_0$ coincides with the apex of the slice and belongs already to the left boundary: the dividing line is therefore
reduced to the backward edge $x_0\to y_0$, hence, in particular, has $p=0$.

For $\ell\geq 3$, we now decompose the $\ell$-slice by cutting it along the dividing line as well as 
along the leftmost backward shortest paths from all the vertices $x_m$ to the apex (taking the backward oriented edge $x_m\to y_m$
as first edge) for $m=1$ to $p-1$ (note that the leftmost backward shortest path which starts from the edge $x_p\to y_p$ sticks 
to the left boundary and needs not being cut). 
This cutting results into $p$ domains in the upper part of the slice and one domain
in the lower part (see figure \ref{fig:dividingone}). The $p$ domains of the upper part are easily recognized as being $\ell_m$-slices with $2\leq \ell_m  \leq \ell-1$ hence,
for arbitrary $\ell\leq k$, each of this $p$ domains is enumerated by $T_{k-1}$. As for the lower part, it is enumerated by some
generating function $\mathcal{K}_p$ which does not depend on $\ell$ and will be discussed below.
We deduce the relation
\begin{equation}
T_k=\mathcal{K}(T_{k-1})\ , \qquad \mathcal{K}(T)\equiv \sum_{p\geq 0} \mathcal{K}_p\, T^p\ ,
\label{eq:relation}
\end{equation} 
where the generating function $\mathcal{K}_p$ of the lower part and the corresponding \emph{kernel} $\mathcal{K}(T)$ will be computed
in the next Section.

As a final remark, let us note that an equivalent way to obtain the above cutting lines 
consists in first transforming our $\ell$-slice into a particular slice with inner faces of degree $4$ 
only, as those considered in \cite{G15b} and then applying on this new slice the decomposition described in \cite{G15b} in this context of
$4$-valent inner faces. As explained above, the transformation from $3$-valent to $4$-valent faces is achieved canonically by simply erasing all the long-edges of the original 
$\ell$-slice. The reader is invited to verify that this alternative way of decomposing our $\ell$-slices matches precisely that of the present approach. 

\subsection{The generating function $\boldsymbol{\mathcal{K}_p}$ and the kernel $\boldsymbol{\mathcal{K}(T)}$}
\label{sec:kernel}    
\begin{figure}
\begin{center}
\includegraphics[width=12cm]{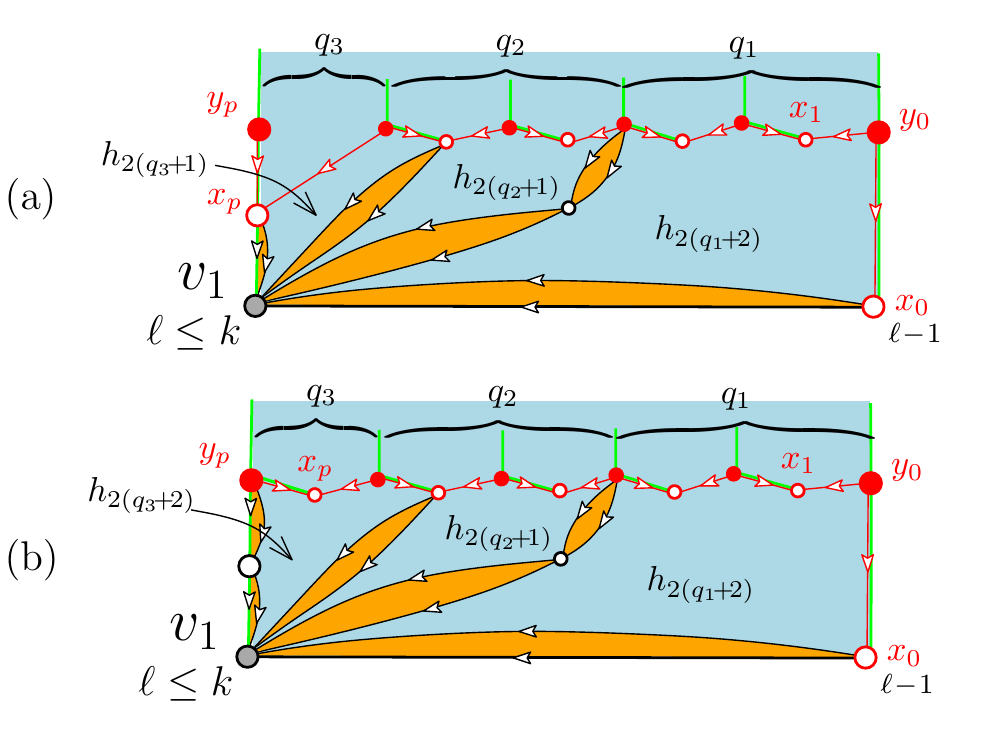}
\end{center}
\caption{A Schematic picture of the decomposition of the part of the slice lying below the dividing line, obtained by marking all
the bundles or pairs of bundles connecting the endpoint $v_1$ of the base to vertices of the dividing line (see text - the light-blue and orange backgrounds are supposedly filled with black and white faces).
Here we have $3$ intermediate domains with $q_1=2$, $q_2=2 $ and $q_3=1$.}
\label{fig:Tknewrec}
\end{figure}
We now come to the generating function $\mathcal{K}_p$ for the lower part of the $\ell$-slice, i.e.\ the domain lying on the same side of the dividing line as
the endpoint $v_1$ of the base. This vertex, at distance $\ell$ from the apex, is connected in general to a number of vertices
$x_{m_1}, x_{m_2},\cdots$ of the dividing line at distance $\ell-1$. This includes in particular the vertex $x_0$ as well
as the vertex $x_p$ when the dividing line precisely hits the left boundary at $x_p$\footnote{A connection from $v_1$ to $x_p$ 
may also be present when the dividing line hits the left boundary at $y_p$ but it is not mandatory in this case.} (see figure \ref{fig:Tknewrec}-(a)).
The connection from $v_1$ to such a vertex $x_{m_j}$ is performed in general by several backward oriented edges and the part of the map
in-between the leftmost and rightmost such edges forms what we called a bundle, as represented in figure \ref{fig:bundle} (recall that
bundles are enumerated by $R_1$). 
Similarly, the vertex $v_1$ may also in general be connected to
a number of vertices
$y_{m_1}, y_{m_2},\cdots$ of the dividing line at distance $\ell-2$ by backward oriented two-step paths passing via some intermediate vertex at distance
$\ell-1$ lying \emph{strictly} in the lower part\footnote{Here we discard possible connections to some $y_m$ via a backward oriented two-step paths having $x_{m-1}$ as intermediate 
vertex.}. This includes in particular the vertex $y_p$ when the dividing line hits the left boundary at $y_p$ (see figure \ref{fig:Tknewrec}-(b)).
The connection to such a vertex $y_{m_j}$ passing through \emph{a fixed intermediate} vertex $z$ is performed in general by several backward oriented edges
from $v_1$ to $z$  and by several backward oriented edges
from $z$ to $y_{m_j}$. The leftmost and rightmost such edges to and from the intermediate vertex delimit a part of the map which now forms a \emph{pair of bundles}. 
Note also that several intermediate vertices $z_1,z_2,\cdots$ may exist for a given $y_{m_j}$, each giving rise to its own pair of bundles.

Cutting out all the bundles from $v_1$ to some $x_{m_j}$ and pairs of bundles from $v_1$ to $y_{m_{j'}}$ leaves us with a number $N\geq 1$ of intermediate domains which are rooted maps
whose outer face has a boundary which is a simple closed curve of length $2i\geq 4$. As displayed in figure \ref{fig:Tknewrec}, this boundary is formed of $v_1$, 
of vertices of the dividing line and of intermediate vertices, connected by successive edges which are either part of the dividing line or lie on the boundary of the
bundles: by convention, the edge towards $v_1$ clockwise around the domain serves as root for the corresponding map. 
A precise characterization of these maps describing the intermediate domains will be given below. At this stage, it is enough to mention that
they are all described by the same internal constraints and are fully characterized by the length $2i$ of their boundary.
We will denote by $h_{2i}\equiv h_{2i}(g)$ their generating function for a fixed boundary length $2i\geq 4$.

To compute $\mathcal{K}_p$, we note that the $p$ domains in the upper part, which are $(k-1)$-slices,  are naturally associated with
the $p$ edges $y_{m-1}\to x_{m}$, $m=1,\cdots,p$, which serve as bases for these slices. 
The $n$-th intermediate domain ($1\leq n\leq N$) has on its boundary a number $q_n$ of such edges and it is easily seen (see figure \ref{fig:Tknewrec} -- see also \cite{G15b} for a more detailed argument) 
that the boundary of this domain is then of length:
\begin{itemize}
\item[-] $2(q_n+1)$ if the left\footnote{By ``left" and ``right", we refer here to the position of the boundary in a representation like in 
figure \ref{fig:Tknewrec}, with the base at the bottom.} boundary of the domain is incident to a single bundle;
\item[-] $2(q_n+2)$ if the left boundary of the domain is incident to a pair of bundles.
\end{itemize}
Note in particular that this length does not depend on whether the right boundary of the domain is incident to a single bundle
or to a pair of bundles. Note also that $q_n$ is necessarily at least $1$ if the left boundary of the domain is incident to a single bundle
but that it can be any non-negative integer if the left boundary of the domain is incident to a pair of bundles\footnote{The
case $q_n=0$ occurs in particular in-between pairs of bundles connecting $v_1$ to the same $y_{m_j}$ via successive intermediate vertices.}.
Summing over all possible configurations of bundles and pairs of bundles, we immediately deduce
\begin{equation*}
\mathcal{K}_p=R_1\, \sum_{N\geq 1} \sum_{q_1,q_2,\cdots, q_N \geq 0 \atop q_1+q_2+\cdots + q_N=p} \prod_{n=1}^{N} \left( R_1\, h_{2(q_n+1)}  +R_1^2 \, h_{2(q_n+2)} \right)
\end{equation*}
where we take the convention $h_2=0$ to suppress the contribution $q_n=0$ in the first term of each factor. Here we decided to attach to the $n$-th intermediate domain 
the weight $R_1$ of the incident single bundle or that, $R_1^2$, of the incident pair of bundles \emph{on its left}.
Note finally the presence of the prefactor $R_1$ corresponding to the bundle linking $v_1$ to $x_0$ which does not lie to the left of any intermediate 
domain.
Summing over $p\geq 0$, we deduce
\begin{equation}
\begin{split}
\mathcal{K}(T)\equiv \sum_{p\geq 0} \mathcal{K}_p T^p& =R_1 \sum_{N\geq 1}\left(\sum_{q\geq 0}\left( R_1 h_{2(q+1)}+R_1^2\, h_{2(q+2)} \right) T^q  \right)^N\\
& = R_1 \sum_{N\geq 1} \left( R_1\left(T+ R_1\right) \, \Phi(T) \right)^N \\
& =\frac{R_1^2 (T+R_1)\, \Phi(T)}{1-R_1(T+R_1)\, \Phi(T)}\ , \qquad \Phi(T)\equiv  \sum_{i\geq 2} h_{2i}\, T^{i-2}\\
\end{split}
\label{eq:Ncal}
\end{equation}
(recall that $h_2=0$). This and \eqref{eq:relation} explain the announced form \eqref{eq:newrecur} of our recursion for $T_k$. 

\subsection{The generating function $\boldsymbol{h_{2i}}$ and the function $\boldsymbol{\Phi(T)}$}
\label{sec:h2i}     
\begin{figure}
\begin{center}
\includegraphics[width=5cm]{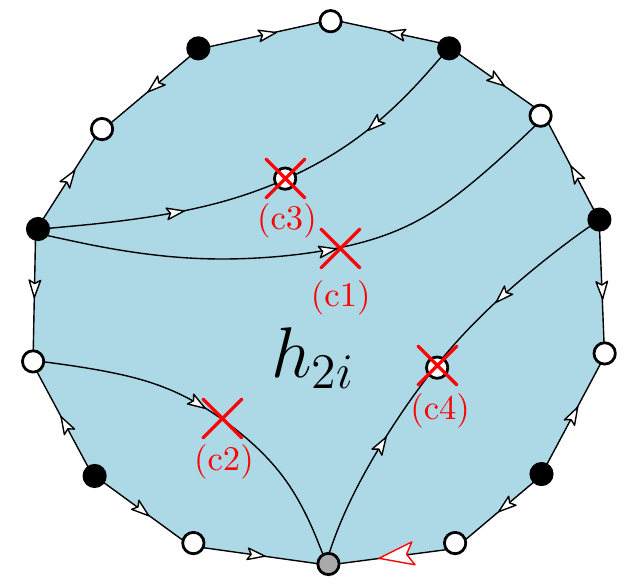}
\end{center}
\caption{A schematic picture of a map enumerated by $h_{2i}$ indicating the edge connections forbidden by the constraints (c1)-(c4) of the text. The length of the boundary of the map is $2i$ (the light-blue background is supposedly filled with black and white faces). }
\label{fig:h2i}
\end{figure}
To compute $\Phi(T)$, we first need to characterize the domains enumerated by $h_{2i}$ for $i\geq 2$. 
As already mentioned, these domains are rooted maps whose outer face has a boundary which is a simple closed curve of length $2i$. 
The inner faces of the maps all have degree $3$ and are colored in black and white so that each edge is incident to a face of each color. The edges are oriented
clockwise around white faces and counterclockwise around black faces. 
To complete our characterization, we first color the vertices in three colors $0$, $1$, $2$ (with $2$ following $1$ following $0$ following $2$ along the oriented 
edges, as explained before), hereafter referred to as \emph{black, white and gray} respectively, in agreement with their representation in the figures. 
The color of the endpoint of the root edge (which corresponds to the vertex $v_1$ in the intermediate domain) in chosen to be \emph{gray}.
Then, as seen from the orientation of edges along the dividing line and at the boundary of the bundles, the colors of the other vertices around 
the boundary of the maps enumerated by $h_{2i}$ alternate between white and black, the two boundary vertices
adjacent to the unique gray boundary vertex being white (in particular, the number of white vertices is $i$ and that of black vertices $i-1$).
Moreover the maps satisfy by construction the following internal constraints, illustrated in figure \ref{fig:h2i}:
\begin{enumerate}
\item[(c1)]{A black vertex of the boundary cannot be linked to a white vertex of the boundary by some oriented edge lying \emph{strictly inside} the map.}
\item[(c2)]{A white vertex of the boundary cannot be linked to the (unique) gray vertex of the boundary by some oriented edge lying strictly inside the map.}
\item[(c3)]{Two distinct black vertices of the boundary cannot have a common white adjacent vertex strictly inside the map.}
\item[(c4)]{The (unique) gray vertex of the boundary cannot have with any black vertex of the boundary a common white adjacent vertex strictly inside the map.}
\end{enumerate}
All these constraint immediately follow from our choice of dividing line and may be verified straightforwardly by the reader.
Let us for instance discuss how constraint (c1) emerges from our construction of the intermediate domains (see figure \ref{fig:Tknewrec}). 
Note that black boundary vertices necessarily come from the dividing line. Assume that an internal edge exists between a black boundary vertex $y_{m_j}$
and a white boundary vertex. Then, if the incident white vertex is also a vertex $x_{m'}$ of the dividing line, either $m'> m$
and a two-step path, made of this internal edge followed by the boundary edge $x_{m'}\to y_{m'}$ lies to the left of the two-step path 
$y_{m}\to x_{m+1}\to y_{m+1}$ originally taken by the dividing line, or $m'\leq m$ and a two-step path, made of the boundary edge $y_{m'-1}\to x_{m'}$ followed the internal edge
lies to the left of the two-step path 
$y_{m'-1}\to x_{m'}\to y_{m'}$ originally taken by the dividing line. In both cases, we have a contradiction.
If instead the white vertex was originally an intermediate vertex in the middle of a pair of bundles on the right side of the domain, the sequence made of 
the boundary edge of the upper bundle followed by the internal edge creates again a two-step path strictly to the left
of the dividing line, a contradiction. If the pair of bundles is on the left side of the domain, an undesirable two-step path
is now obtained by taking the internal edge followed by the boundary edge of the upper bundle. 

\begin{figure}
\begin{center}
\includegraphics[width=5cm]{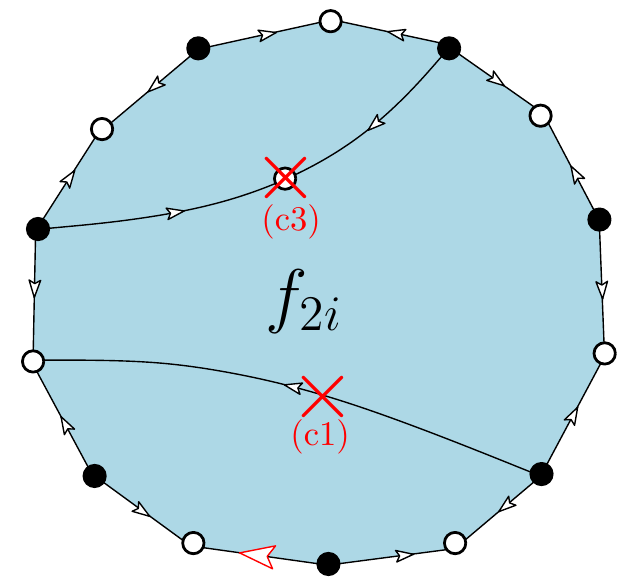}
\end{center}
\caption{A schematic picture of a map enumerated by $f_{2i}$ indicating the edge connections forbidden by the constraints (c1) and (c3) of the text. The length of the boundary of the map is $2i$ (the light-blue background is supposedly filled with black and white faces).}
\label{fig:f2i}
\end{figure}
Having defined the maps enumerated by $h_{2i}$, we may now determine $\Phi(T)=\sum_{i\geq 2}h_{2i} T^{i-2}$
by writing some self-consistency relation for this quantity. In order to get a simple closed relation for $\Phi(T)$, we have to introduce 
yet another family of rooted maps whose outer face has again a boundary made of a simple closed curve of length $2i$, $i\geq 2$, 
but whose boundary color assignment is  slightly different. 
The generating functions of these maps will be denoted by $f_{2i}\equiv f_{2i}(g)$, $i\geq 2$ and are gathered in the quantity
\begin{equation}
\Omega(T)\equiv \Omega(T,g)=\sum_{i\geq 2} f_{2i}(g)\, T^{i-2}\ .
\label{eq:Omegadef}
\end{equation} 
In maps enumerated by $f_{2i}$, the endpoint of the root edge is white and the colors around the boundary 
alternate between white and black (see figure \ref{fig:f2i} -- 
in particular, there are $i$ white boundary vertices and $i$ black boundary vertices).
The maps are moreover required to \emph{satisfy the same constraints} (c1) \emph{and} (c3) \emph{above}, as illustrated in figure \ref{fig:f2i}
(the constraints (c2) and (c4) are pointless since there is no gray boundary vertex).
\vskip 1.9cm
\section{Computation of $\Phi(T)$ and solution of the recursion}
\label{sec:phi}  
\subsection{A closed system for $\boldsymbol{\Phi(T)}$ and $\boldsymbol{\Omega(T)}$}
\label{sec:system}
\begin{figure}
\begin{center}
\includegraphics[width=11cm]{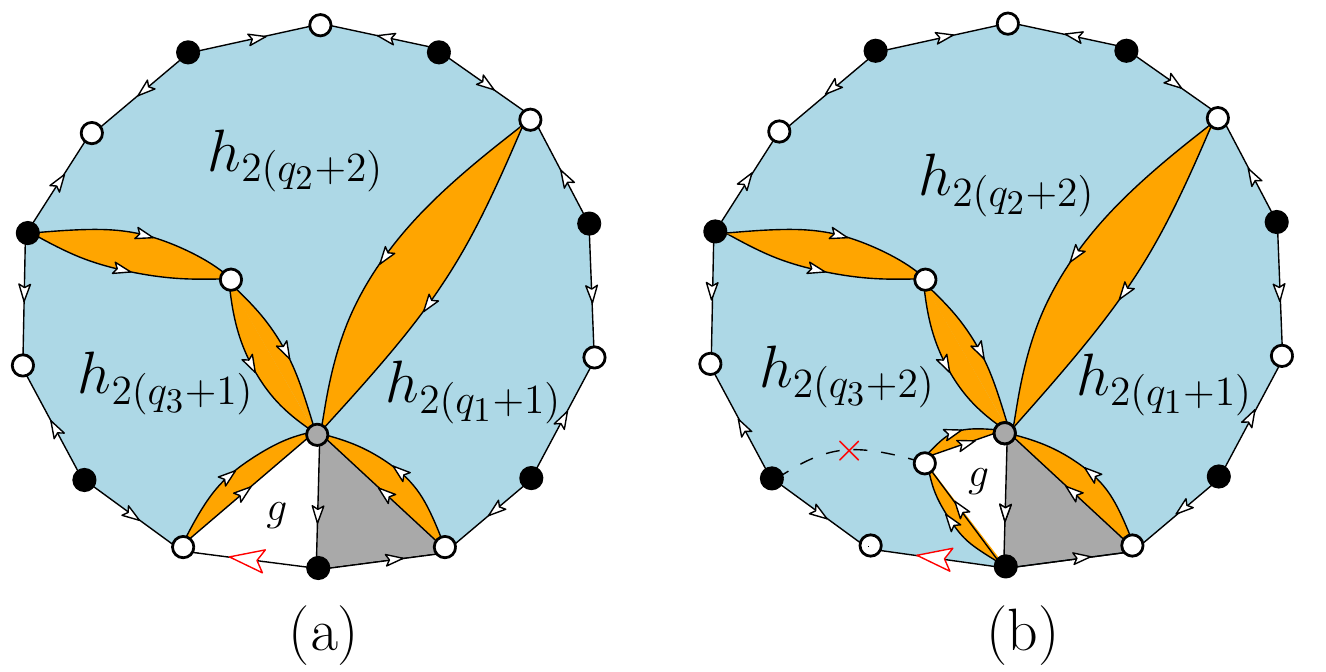}
\end{center}
\caption{Decomposition of a map enumerated by $f_{2i}$ (see tex -- the light-blue and orange backgrounds are supposedly filled with black and white faces). Here the $q_m$'s  
designate the number of black to white counterclockwise oriented boundary edges incident to the domains at hand.}
\label{fig:Omega}
\end{figure}
The functions $\Phi(T)$ and $\Omega(T)$ may be determined in terms of $R_1$ and $g$ by the following closed system: 
\begin{equation}
\begin{split}
\Omega(T)& = \frac{g}{T}\left\{\frac{R_1^2(T+R_1)\Phi(T)}{1-R_1(T+R_1)\Phi(T)}-\frac{R_1^3\,  h_4}{1-R_1^2\,  h_4}\right\}\\
\Phi(T)& =g +\frac{g}{T}\left\{\frac{R_1^2\, T\Omega(T)}{1-R_1(T+R_1)\Omega(T)}+\frac{R_1^3}{1-R_1^2\, h_4}\bigg(\frac{(\Phi(T)-h_4)}{1-R_1(T+R_1)\Omega(T)}\right.
\\ & \hskip 8.cm
\left.+ \frac{R_1\, T \Omega(T)\, h_4}{1-R_1(T+R_1)\Omega(T)}\bigg)
\right\} \ .\\
\end{split}
\label{eq:PhiOmega}
\end{equation}
Note that this system also fixes the value of $h_4=\Phi(0)$ as a function of $g$.

The first equation in this system is a direct consequence of the identification
\begin{equation}
f_{2i}=g\, \mathcal{K}_{i-1}\ , \qquad i\geq 2
\label{eq:fiN}
\end{equation}
which may be understood as follows (see figure \ref{fig:Omega} for an illustration): consider, in a map
enumerated by $f_{2i}$ the boundary edge starting from the (black) origin of the root-edge \emph{counterclockwise} around the map.
It leads to a white boundary vertex (necessarily different from the endpoint of the root edge since the length of the boundary is 
$2i\geq 4$) and has a black inner face on its left whose third vertex is gray an therefore cannot lie on the boundary. The face on the
right
of the edge incident to this black inner face and leading from the gray vertex back to the black origin of the root-edge is a white inner face whose third vertex is white and may
either be the endpoint of the root edge (case (a) in figure \ref{fig:Omega}) or another white vertex (case (b) in figure \ref{fig:Omega}).
In this latter case, since this white vertex is connected to the black origin of the root-edge, it cannot, because of
(c1), belong to the boundary and it cannot, because of (c3), be connected to another black vertex on the boundary.
If we now \emph{remove the black and white inner faces} that we just singled out (which contribute a weight $g$ to $f_{2i}$), the rest of the map (with 
boundary length $2i$ in case (a) and $(2i+1)$ in case (b)) 
is characterized by exactly the same constraints as those defining $\mathcal{K}_p$ and illustrated in figure \ref{fig:Tknewrec}
(with the (a) and (b) cases there matching precisely the present (a) and (b) cases)  
with $p=i-1$ (since, as seen in figure \ref{fig:Tknewrec}, the boundary of the lower part of the slice has length $2p+2$ in case
(a) and $2p+3$ is case (b)).
This identification is best seen by drawing the bundles performing the direct connections from the singled out gray vertex to white boundary vertices
as well as the pairs of bundles performing the connections from this singled out gray vertex to black boundary vertices
via some intermediate white bulk vertex. We recover as in figure \ref{fig:Tknewrec} a number $N$ of intermediate domains enumerated by 
$h_{2q_m+1}$ or $h_{2q_m+2}$ according to whether the domain is followed counterclockwise by a single bundle or by a pair of bundles
respectively ($q_m$ being the number of black to white counterclockwise oriented boundary edges incident to the domain at hand).

 Plugging \eqref{eq:fiN} into \eqref{eq:Omegadef} allows us to write
\begin{equation*}
\Omega(T)=\frac{g}{T}(\mathcal{K}(T)-\mathcal{K}(0))
\end{equation*}
and, using the explicit form \eqref{eq:Ncal} of $\mathcal{K}(T)$, we immediately arrive at the desired equation \eqref{eq:PhiOmega} for $\Omega(T)$.

\begin{figure}
\begin{center}
\includegraphics[width=11cm]{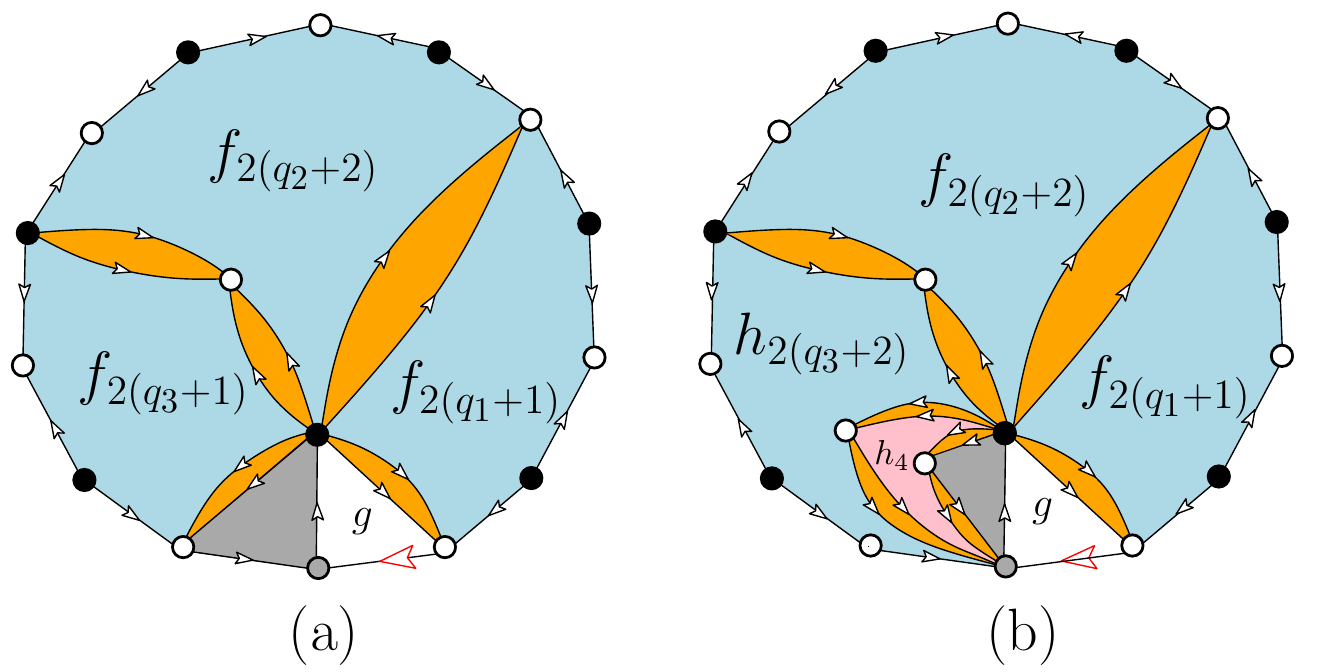}
\end{center}
\caption{Decomposition of a map enumerated by $h_{2i}$ (see text). In case (b), the last $M$ intermediate domains (in pink) correspond to domains enumerated by $h_4$ and delimited on both sides by pairs of bundles
connecting the singled out black vertex in the bulk to the gray boundary vertex (here $M=1$ --  the light-blue, pink and orange backgrounds are supposedly filled with black and white faces).}
\label{fig:Phi}
\end{figure}

To obtain the second equation in \eqref{eq:PhiOmega}, we now consider a map
enumerated by $h_{2i}$. Its root edge now points from a white boundary vertex to the unique gray boundary vertex and has a white inner face on its right 
whose third vertex is black. The edge incident to this white inner face and leading from the gray boundary vertex to the singled out black vertex cannot be a 
boundary edge (there are no gray $\to$ black boundary edges), hence it has on its left a black inner face whose third vertex is white. 
This white vertex may either be the boundary vertex preceding the gray boundary vertex counterclockwise around the boundary (case (a) in figure \ref{fig:Phi}) or be 
another white vertex (case (b) in figure \ref{fig:Phi}).  
In this latter case, since this white vertex is connected to the gray boundary vertex, it cannot, because of (c2), belong to the boundary and it cannot,
because of (c4), be connected to a black boundary vertex.

\begin{figure}
\begin{center}
\includegraphics[width=4cm]{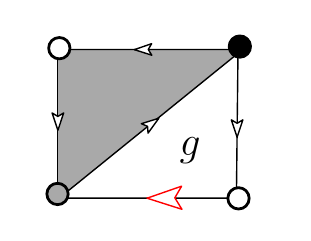}
\end{center}
\caption{The simplest map contributing to $h_4$. The contribution of this map to $\Phi(T)$ is $g$. }
\label{fig:single}
\end{figure}

Note that the singled out black vertex may itself be a boundary vertex but, because of constraints (c1), (c2) and (c4), the map necessarily reduces in this case to 
the trivial map of boundary length $4$ displayed in figure \ref{fig:single}. This map contributes a weight $g$ to $h_4$, hence to $\Phi(T)$ and this explains 
the first term in the second line of \eqref{eq:PhiOmega}.
In all other cases, the singled out black vertex does not lie on the boundary. As before, we may then remove the singled out black and white faces (which contribute a weight $g$ to $h_{2i}$) 
and draw the bundles performing direct connections from the singled out black vertex to white boundary vertices
as well as the pairs of bundles performing connections from the singled out black vertex to black boundary vertices
via some intermediate white bulk vertex (see figure \ref{fig:Phi}). This defines a number $N\geq 1$ of intermediate domains. In case (a), all these domains are
enumerated by $f_{2(q_m+1)}$ or $f_{2(q_m+2)}$ according to whether the domain is followed counterclockwise by a single bundle or by a pair of bundles
respectively ($q_m$ begin the number of black to white counterclockwise oriented boundary edges incident to the domain at hand), \emph{the $N$-th domain 
being necessarily followed by a single bundle}. Repeating the argument that led to \eqref{eq:Ncal}, the total contribution to $\Phi(T)$ of these situations is\footnote{\label{foot1}The
prefactor $g/T$ includes a weight $g$ for the removed pair of inner faces and a factor $1/T$ to compensate the fact that the required power $i-2$ of $T$ in $\Phi(T)$ is $1$ 
unit less than the sum of the $q_m$'s (i.e.\ $\sum_{m=1}^N q_m=i-1$).}
\begin{equation*}
\begin{split}
&\frac{g}{T}\left\{R_1 \sum_{N\geq 1}\left(\sum_{q\geq 0}\left( R_1 f_{2(q+1)}+R_1^2\, f_{2(q+2)} \right) T^q  \right)^{N-1}\times \left(\sum_{q\geq 0}R_1 f_{2(q+1)}\,  T^q  \right)\right\}\\
=&\frac{g}{T}\left\{  R_1 \sum_{N\geq 1} \left( R_1\left(T+ R_1\right) \, \Omega(T) \right)^{N-1} \times R_1\, T\Omega(T) \right\}\\
=& \frac{g}{T}\left\{\frac{R_1^2\, T\Omega(T)}{1-R_1(T+R_1)\Omega(T)}\right\}\\
\end{split}
\end{equation*}
(with the convention $f_2=0$). This explains the second term in the second line of \eqref{eq:PhiOmega}.

In case (b), among the $N$ domains, the last $M$ ones ($0\leq M\leq N-1$) corresponds to domains delimited on both sides by pairs of bundles
connecting the singled out black vertex to the gray boundary vertex and contribute a weight $R_1^2 h_4$ each (see figure \ref{fig:Phi}-(b)). Altogether, these domains
eventually contribute (after summation over $M$) a prefactor 
\begin{equation*}
\frac{1}{1-R_1^2\, h_4}
\end{equation*}
to the desired contribution of case (b). As for the remaining $N-M$ domains, 
they are enumerated 
enumerated by $f_{2(q_m+1)}$ or $f_{2(q_m+2)}$ according to whether the domain is followed counterclockwise by a single bundle or a pair of bundles
respectively ($q_m$ begin the number of black to white counterclockwise oriented boundary edges incident to the domain at hand), \emph{except for 
the $(N-M)$-th domain} which is necessarily followed by a pair of bundles and is enumerated by $h_{2(q_{N-M}+2)}$.  
\begin{figure}
\begin{center}
\includegraphics[width=14cm]{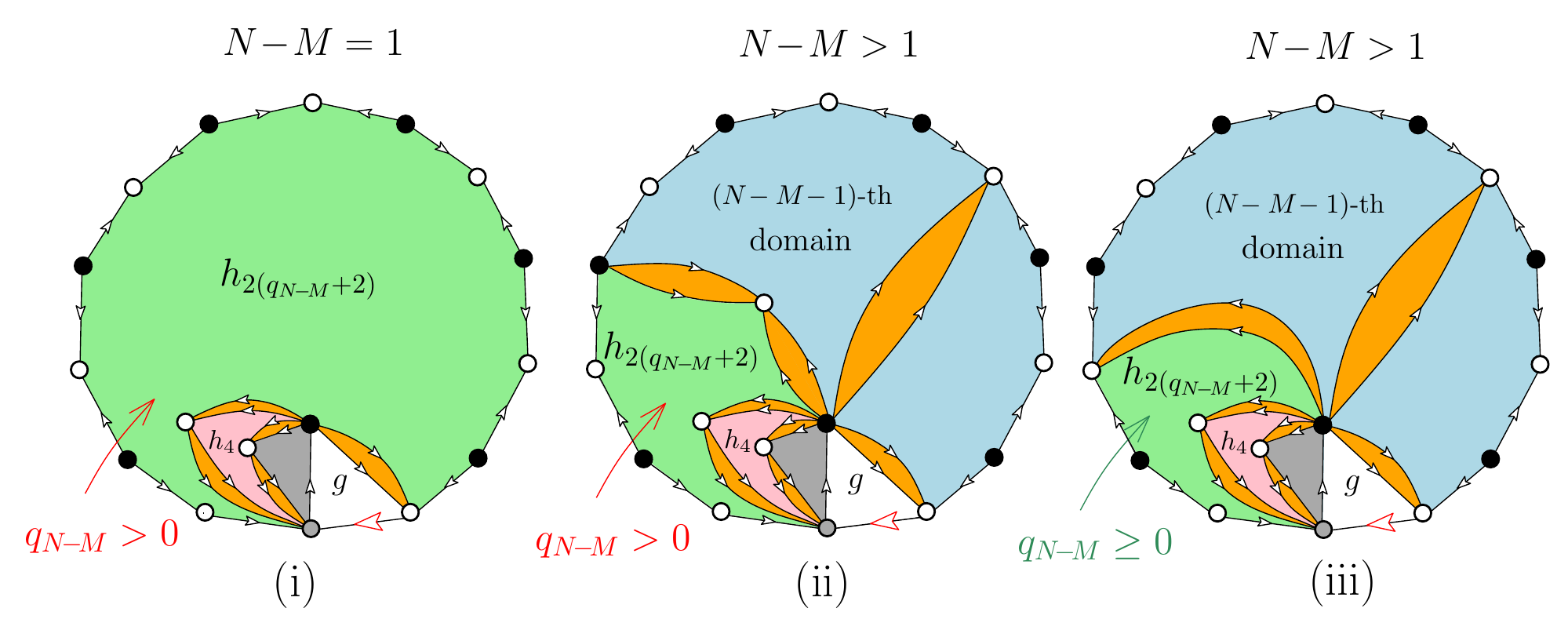}
\end{center}
\caption{Particular classes of maps enumerated by $h_{2i}$ giving rise to constraints on the boundary length of the $(N-M)$-th intermediate domain, here in green 
(see text for details -- as in figure \ref{fig:Phi}, the last $M$ domains, in pink, correspond to domains delimited on both sides by pairs of bundles
connecting the singled out black vertex in the bulk to the gray boundary vertex).  In case (i) and (ii) the number $q_{N-M}$ of black to white counterclockwise oriented boundary edges incident to the domain at hand must satisfy $q_{N-M}>0$ while in case (iii), $q_{N-M}=0$ is also allowed (the light-blue, light-green, pink and orange backgrounds are supposedly filled with black and white faces).}
\label{fig:Philengthfour}
\end{figure}
Repeating the summation argument above (summing now over $(N-M)$), these domains contribute naively a total weight\textsuperscript{\ref{foot1}}
\begin{equation*}
\frac{g}{T}\left\{\frac{R_1^3\, \Phi(T)}{1-R_1(T+R_1)\Omega(T)}\right\}
\end{equation*}
to $\Phi(T)$. This expression is however not fully satisfactory since it assumes that the $(N-M)$-th domain may have an arbitrary boundary length 
$2(q_{M-N}+2)\geq 4$ while a length $4$ (i.e.\ $q_{M-N}=0$) is in fact \emph{not allowed} in the following situations: whenever (i) $(N-M)=1$ (i.e.\ the domains enumerated by $f_{2(q_m+1)}$ or $f_{2(q_m+2)}$ are absent) since
this would lead to an original map with boundary length $2i=2$, not present in $\Phi(T)$ or whenever (ii) $(N-M)> 1$ and the $(N-M-1)$-th domain
is followed counterclockwise by a pair of bundles, since, by definition of $M$, this pair of bundles cannot hit the boundary at the gray vertex (see figure \ref{fig:Philengthfour}
for illustration). On the contrary, the 
length $2(q_{M-N}+2)=4$ (i.e.\ $q_{M-N}=0$) is allowed whenever (iii) $(N-M)> 1$ and the $(N-M-1)$-th domain
is followed counterclockwise by a single bundle, since this bundle may hit the boundary at the white vertex preceding the gray vertex counterclockwise
around the map. Taking these corrections into account, our naive estimate but be modified into
\begin{equation*}
\frac{g}{T}\left\{\frac{R_1^3\, (\Phi(T)-h_4)}{1-R_1(T+R_1)\Omega(T)}+ \frac{R_1^3\, h_4\times R_1T \Omega(T)}{1-R_1(T+R_1)\Omega(T)}\right\}
\end{equation*}
where the first term simply forbids the length $2(q_{M-N}+2)=4$ in all cases (i), (ii) and (iii) while the second term reintroduces it in the allowed case (iii) .
Incorporating the prefactor $1/(1-R_1^2\, h_4)$ above, we eventually obtain the third and fourth terms in the second line of \eqref{eq:PhiOmega}.

\subsection{Solution of the system (\ref{eq:PhiOmega})}
\label{sec:solsystem}
We may slightly simplify the system \eqref{eq:PhiOmega} by performing the rescalings\footnote{\label{foot2}Here we view this rescaling as a simple change of variables but the reader may easily verify,
by a simple substitution at the level of maps, that $\tilde{h}_{2i}(G)$ (respectively $\tilde{f}_{2i}(G)$) actually enumerates the sub-family
of the maps enumerated by $h_{2i}(g)$ (respectively $f_{2i}(g)$) made of those maps having no multiple edges
connecting a black to a white vertex or a white to a gray vertex (multiple edges connecting a gray to a black vertex are allowed),
with now a weight $G$ per white face.}
\begin{equation*}
\begin{split}
& G\equiv g\, R_1^2\ , \qquad  t\equiv \frac{T}{R_1}\ , \qquad \tilde{h}_{2i}(G)\equiv R_1^{i}\, h_{2i}(g)\ , \qquad \tilde{f}_{2i}(G)\equiv R_1^{i}\, f_{2i}(g)\ , \\
& 
\hskip -1.2cm \phi(t)\equiv \phi(t,G)=\sum_{i\geq 2} \tilde{h}_{2i}(G)\, t^{i-2}
= R_1^{2}\, \Phi(T)
\ , \quad \omega(t)\equiv \omega(t,G)=\sum_{i\geq 2} \tilde{f}_{2i}(G)\, t^{i-2}
= R_1^{2}\, \Omega(T)\ .\\
\end{split}
\end{equation*}
With these new variables, our system simplifies into
\begin{equation}
\begin{split}
\omega(t)& = \frac{G}{t}\left\{\frac{(t+1)\phi(t)}{1-(t+1)\phi(t)}-\frac{\tilde{h}_4}{1-\tilde{h}_4}\right\}\\
\phi(t)& =G +\frac{G}{t}\left\{\frac{t\omega(t)}{1-(t+1)\omega(t)}+\frac{1}{1-\tilde{h}_4}\left\{\frac{\phi(t)-\tilde{h}_4+t\, \omega(t)\, \tilde{h}_4}{1-(t+1)\omega(t)}\right\}
\right\} \\
\end{split}
\label{eq:phiomega}
\end{equation}
where the quantity $R_1$ dropped out. Both $\phi(t)$ and $\omega(t)$ implicitly depend on $G$, and so does $\tilde{h}_4$.

As for the recursion relation \eqref{eq:newrecur} itself, it also gets simpler by using similar rescalings:
\begin{equation*}
r_k\equiv \frac{R_k}{R_1}\ , \qquad  t_k\equiv \frac{T_k}{R_1}=r_k -1\ .
\end{equation*}
With these variables, \eqref{eq:newrecur} translates into
\begin{equation}
t_k=\frac{(t_{k-1}+1)\, \phi(t_{k-1})}{1-(t_{k-1}+1)\, \phi(t_{k-1})}\ .
\label{eq:newrecursimp}
\end{equation}

Solving a system of the form \eqref{eq:phiomega} a standard exercise. The first step consists in getting an expression for $\tilde{h}_4$ in terms
of $G$ as follows: from the first and second equations in \eqref{eq:phiomega}, we immediately deduce, for a given $G$, the two relations:
\begin{equation}
\begin{split}
\tilde{h}_4& =\tilde{h}_4^{(1)}(t,\phi(t),\omega(t))\ \ \hbox{with}\ \ \tilde{h}_4^{(1)}(t,\phi,\omega)= \frac{G (t+1)\, \phi +t\, \omega\, ((t+1) \phi -1)}{G+t\, \omega\, ((t+1) \phi -1)}\ ,\\
\tilde{h}_4& =\tilde{h}_4^{(2)}(t,\phi(t),\omega(t))\ \ \hbox{with}\ \ \tilde{h}_4^{(2)}(t,\phi,\omega)= \frac{G \left(\omega  t^2-t-\phi \right)-t\, \phi\,  ((t +1)\omega -1)}{G (t+1) (t\, \omega -1)-t\, \phi \, ((t +1)\omega  -1)}
\ .\\
\end{split}
\label{eq:h4exprl}
\end{equation} 
Since $\tilde{h}_4$ does not depend on $t$, we may then write
\begin{equation*}
\begin{split}
0& =\frac{d\tilde{h}_4}{dt}=\frac{\partial \tilde{h}_4^{(1)}}{\partial \phi}\phi'(t)+\frac{\partial \tilde{h}_4^{(1)}}{\partial \omega}\omega'(t)+\frac{\partial \tilde{h}_4^{(1)}}{\partial t}\\
0& =\frac{d\tilde{h}_4}{dt}=\frac{\partial \tilde{h}_4^{(2)}}{\partial \phi}\phi'(t)+\frac{\partial \tilde{h}_4^{(2)}}{\partial \omega}\omega'(t)+\frac{\partial \tilde{h}_4^{(2)}}{\partial t}\\
\end{split}
\end{equation*} 
and, upon eliminating $\omega'(t)$,
\begin{equation}
0=\left\{\frac{\partial \tilde{h}_4^{(1)}}{\partial \phi}\frac{\partial \tilde{h}_4^{(2)}}{\partial \omega} -\frac{\partial \tilde{h}_4^{(2)}}{\partial \phi}\frac{\partial \tilde{h}_4^{(1)}}{\partial \omega}  \right\}\phi'(t)+ 
\left\{ \frac{\partial \tilde{h}_4^{(1)}}{\partial t}\frac{\partial \tilde{h}_4^{(2)}}{\partial \omega} -\frac{\partial \tilde{h}_4^{(2)}}{\partial t}\frac{\partial \tilde{h}_4^{(1)}}{\partial \omega} \right\}\ .
\label{eq:frac}
\end{equation}
To obtain $\tilde{h}_4$ in terms of $G$, we may now solve the system \eqref{eq:phiomega} on some particular line $t=t(G)$ where the second term between braces in \eqref{eq:frac} vanishes. 
From the explicit expressions \eqref{eq:h4exprl}, it is easily checked that such vanishing occurs whenever $\omega\equiv\omega(t(G))$ is related to $\phi\equiv \phi(t(G))$ via
\begin{equation}
\omega=  \frac{t(t+1)\, \phi+G\, (t+1)-t}{t(t+1)((t+1)\, \phi  -1)}\ .
\label{eq:omegaval}
\end{equation}
On the line $t=t(G)$, the first term between braces in \eqref{eq:frac} must also vanish (since, generically, $\phi'(t(G))\neq 0$). Plugging
the above value \eqref{eq:omegaval} for $\omega$, this leads to the following equation for $\phi$:
\begin{equation*}
\hskip -1.cm 0=(t+1) G^2+\left(2 (t+1)^2 \phi ^2-3 (t+1) \phi +1\right) G-\phi  \left((t+1)\, \phi -1\right)^2\ ,
\end{equation*}
which we complete by demanding that $\tilde{h}_4^{(1)}=\tilde{h}_4^{(2)}$, namely, after plugging again the value \eqref{eq:omegaval} of $\omega$:
\begin{equation*}
\begin{split}
&\hskip -1.cm 0= (t+1)^2 \left((t+1)^2\, \phi  +1-t \right) G^2 -(t+1) ((t +1)\,\phi-1) \left(2 t(t+1)\, \phi -2 t+1\right) G\\
& \hskip 10.cm + t ((t+1)\, \phi-1)^3 . \\
\end{split}
\end{equation*}
Imposing the two equations above determines an equation for the line $t=t(G)$ as well as the value of $\phi$ on this line. The value of $\omega$
is then obtained through \eqref{eq:omegaval}. After some straightforward calculations, we find explicitly:
\begin{equation}
0=(t+1)^5\, G^2 -(t-1) (t+1)^2\, G -t
\label{eq:tG}
\end{equation}
and the values
\begin{equation*}
\begin{split}
& \phi= \frac{t}{(t+1)^2}\ ,
\\&
\omega= \frac{t-G (t+1)^2}{t (t+1)}\ .
\\
\end{split}
\end{equation*}
The reader should remain aware that these latter expressions for $\phi$ and $\omega$ do not hold for arbitrary $t$ and $G$ but \emph{only on the line $t=t(G)$} implicitly defined
by \eqref{eq:tG}. Still, upon plugging these values into the expression \eqref{eq:h4exprl} for $\tilde{h}_4^{(1)}$ (or equivalently for $\tilde{h}_4^{(2)}$ since we imposed
$\tilde{h}_4^{(1)}=\tilde{h}_4^{(2)}$),
we obtain the desired value of $\tilde{h}_4$, namely:
\begin{equation*}
\tilde{h}_4= \frac{G^2 t (t+1)^3+t-G \left(t^3+3 t^2+2 t+1\right)}{G^2 (t+1)^4-G (t+1)^3+t}\ ,
\end{equation*}
valid on the line $t=t(G)$. 

Here, $t$ may be viewed as a simple parametrization of $G$ via \eqref{eq:tG}. In order to have a slightly simpler parametric expression for $\tilde{h}_4(G)$, 
and in view of the explicit form (quadratic in $G$) of the relation \eqref{eq:tG}, 
we decide to parametrize $t$ itself it by the quantity $C$
solution of the quadratic equation $0=(t+1)\, C^2-(t-1)\, C-t$. In other words, we set
\begin{equation*}
t=-\frac{C (C+1)}{C^2-C-1}
\end{equation*} 
so that $G$ and $\tilde{h}_4(G)$ are now expressed as rational functions in $C$, namely:
\begin{equation}
\begin{split}
& \tilde{h}_4(G)= \frac{C \left(C^3+2 C^2-C-1\right)}{(C-1) (2 C+1)^2} \\
&\hbox{where}\ \ G= \frac{C \left(C^2-C-1\right)^2}{(2 C+1)^2} \ .\\
\end{split}
\label{eq:h4Gparam}
\end{equation}
The above parametrization is univocal for $0\leq C\leq 1/2$, leading to a value of $G$ in the range 
$0\leq G\leq 25/128$. We find in particular from this parametrization the small $G$ expansion
\begin{equation*}
\tilde{h}_4(G)=G+G^3+3 G^4+9 G^5+31 G^6+114 G^7+435 G^8+1713 G^9+6924 G^{10}+O\left(G^{11}\right)\ .
\end{equation*}

Going back to the system \eqref{eq:phiomega}, now \emph{for arbitrary $t$ and $G$}, we may plug the 
above parametrization of $G$ and $\tilde{h}_4$ to obtain, after eliminating $\omega=\omega(t)$, an equation for
$\phi=\phi(t)$ as a function of $t$ and the quantity $C$ (parametrizing $G$ as above), namely:
\begin{equation*}
\begin{split}
\hskip -1.cm 0= &\  C\, (C+1)^3 \bigg(C\, (C^3+2C^2-C-1)+(C^2-C-1)^2\, t \bigg)\\
&  -(2 C+1)^2 \bigg(C\, (C-1) (C+1)^3+(C+1) (2 C^4-2 C^3+C^2+3 C+1)\, t\\
& \hskip 4.cm + C \Big(C^2-C-1\Big)^2\, t^2 \bigg)\, \phi +(2 C+1)^4 \, t\, (t+1)\, \phi^2 \ .\\
   \end{split}
\end{equation*}
This is a quadratic equation in $\phi$, whose discriminant factorizes into 
\begin{equation*}
\begin{split}
&\Delta  =(2 C+1)^4 \left(C\, (C+1)+(C^2-C-1)\, t \right)^2 \, \times \delta \ ,\\
& \delta=\left((C-1)^2 (C+1)+C \left(C^2-C-1\right)\, t \right)
\left((C+1)^3+C \left(C^2-C-1\right)\, t \right)\ . \\
       \end{split}
   \end{equation*}
We now note that $\delta$ is itself the discriminant of the following quadratic equation\footnote{\label{foot3}This choice of the parameter $Y$ is 
definitely not unique and many other possibilities exist, based on other choices for the quadratic 
equation (with still coefficients linear in $t$) having discriminant $\delta$. The various
steps of our calculation may easily be repeated with these alternative definitions of $Y$ and should eventually lead to the same final expression for $\phi(t)$.} in the variable $Y$, with coefficients
linear in $t$ :
\begin{equation}
0= C^2 (C+1)^2+\left((C+1) \left(C^2+1\right)+C \left(C^2-C-1\right) t\right) Y +Y^2\ . 
\label{eq:eqY}
\end{equation}
This suggest to parametrize $t$ by this quantity $Y$, namely set
\begin{equation}
t=-\frac{(C+Y+1) \left(C^3+C^2+Y\right)}{C \left(C^2-C-1\right) Y}\ .
\label{eq:tY}
\end{equation}
With this parametrization, $\phi$ is now a rational function of $Y$, namely\footnote{The equations for $\phi$ and that for $Y$ are quadratic 
hence have two solutions $\phi_1$, $\phi_2$ and $Y_1$, $Y_2$ respectively. The relation that we give here corresponds to 
one choice of passage from $Y$ to $\phi$, say from $Y_1$
to $\phi_1$ and from $Y_2$ to $\phi_2$. The passage from $Y_1$ to $\phi_2$ and from $Y_2$ to $\phi_1$ is achieved by some
easily computable different rational function, which we could chose as well in our calculation. This would then require picking 
the other branch of $Y$ in the subsequent argument.} 
\begin{equation}
\phi=-\frac{C \left(C^2-C-1\right) Y \left(C^4+2 C^3-Y C^2+C^2+Y C+Y\right)}{(2 C+1)^2 \left(C^2+Y\right) \left(C^3+C^2+Y\right)}\ .
\label{eq:phiY}
\end{equation}
Among the two solutions $Y$ of the quadratic equation \eqref{eq:eqY}, we must pick the one which yields $\phi=\tilde{h}_4$ 
for $t=0$, where $\tilde{h}_4$ is given by \eqref{eq:h4Gparam}. This yields the determination
\begin{equation}
\begin{split}
Y(t)=& -\frac{1}{2} \bigg((C+1) \left(C^2+1\right)+C \left(C^2-C-1\right) t\\& \hskip 2.cm +\sqrt{\left((C+1) \left(C^2+1\right)+C
   \left(C^2-C-1\right) t\right)^2-4 C^2 (C+1)^2}\bigg)\\
   \end{split}
   \label{eq:Yvalt}
 \end{equation}
which, plugged into \eqref{eq:phiY}, yields an explicit expression for $\phi(t)$. In practice, this expression is not really needed 
and the parametrizations \eqref{eq:tY} and \eqref{eq:phiY} are sufficient.
 
 \subsection{Solution of the recursion relation}
 \label{sec:solrecur}
In order to incorporate the above expressions in our recursion relation, we introduce the quantities
\begin{equation*}
Y_k\equiv Y(t_k)\ , \quad \phi_k\equiv \phi(t_k)\ ,
\end{equation*}
which, according to \eqref{eq:tY} and \eqref{eq:phiY}, are linked by the relations
\begin{equation*}
\begin{split}
& t_k=-\frac{(C+Y_k+1) \left(C^3+C^2+Y_k\right)}{C \left(C^2-C-1\right) Y_k}\ ,\\
& \phi_k=-\frac{C \left(C^2-C-1\right) Y_k \left(C^4+2 C^3-Y_k C^2+C^2+Y_k C+Y_k\right)}{(2 C+1)^2 \left(C^2+Y_k\right) \left(C^3+C^2+Y_k\right)}\ .\\
\end{split}
\end{equation*}
Plugging this parametrization of $t_k$ and the corresponding value of $\phi_k$ in our recursion relation \eqref{eq:newrecursimp} 
transforms it into a recursion relation for $Y_k$ itself, which 
takes the remarkably simple factorized form
\begin{equation*}
\begin{split}
0= & \left(C^3(C+1)-C^2\, Y_{k-1}-C^2\, Y_k -Y_{k-1} Y_k\right) \times \\
& \qquad \left(C^2(C+1)^2+(C+1)^2 Y_{k-1} -C(C+1)Y_k + Y_{k-1}
   Y_k\right)\ .\\
   \end{split}
   \end{equation*}
To decide which factor actually vanishes for the correct determination \eqref{eq:Yvalt} of $Y_k=Y(t_k)$, we note that, for $G\to 0$, we have $t_k=O(G)$ and $C=O(G)$ hence, from \eqref{eq:Yvalt}, $Y_k=-1+O(G)$ and 
$Y_{k-1}=-1+O(G)$ as well. This selects the second factor and our recursion relation translates eventually into the following \emph{homographic recursion}:
\begin{equation*}
Y_k=(C+1)^2
\frac{Y_{k-1}+C^2}{C(C+1)-Y_{k-1}}
\end{equation*}
with, from \eqref{eq:Yvalt}, initial condition $Y_1=-(C+1)$ since $t_1=0$.

We now recall how to solve a general homographic recursion of the form:
\begin{equation*}
Y_k=f(Y_{k-1})\ , \qquad f(Y)\equiv \frac{a\, Y+b}{c\, Y+d}\ .
\end{equation*}
This requires introducing the two fixed points $\alpha$ and $\beta$ of the function $f$ (solutions of $f(Y)=Y$) and
the quantity 
\begin{equation*}
W_k=\frac{Y_k-\alpha}{Y_k-\beta}\ .
\end{equation*}
It is indeed easily seen that $W_k$ satisfies a geometric recursion relation $W_k=x \, W_{k-1}$
with 
\begin{equation*}
x \equiv \frac{c\, \beta+d}{c\, \alpha+d}\ .
\end{equation*}
This yields $W_k=x^{k-1}\, W_1$ and by a simple inversion (assuming $\alpha\neq \beta$, hence $x\neq 1$)
\begin{equation*}
Y_k=\frac{\alpha-\beta\, x^{k-1}\, W_1}{1-x^{k-1}\, W_1}=\frac{Y_1  (\beta \, x^{k-1}-\alpha) - \alpha  \beta \,  (x^{k-1}-1)}{Y_1 (x^{k-1}-1)-(\alpha\,  x^{k-1}-\beta)
}\ .  
\end{equation*}
In the present case, we may take
\begin{equation*}
a=(C+1)^2\ , \qquad b= C^2(C+1)^2\ , \qquad c=-1\ , \qquad d=C(C+1)\ ,
\end{equation*}
so that we find
\begin{equation*}
\begin{split}
\alpha & =  \frac{1}{2}(1+C)\left(-\sqrt{1-4 C^2}-1\right)\ ,\\
\beta & =\frac{1}{2}(1+C)\left(\sqrt{1-4 C^2}-1\right)\ ,\\
x & = \frac{1-\sqrt{1-4 C^2}}{2C}\ ,\\
\end{split}
\end{equation*}
with in particular $0\leq x < 1$ for $0\leq C< 1/2$ (i.e.\ $0\leq G< 25/128$).
Upon inverting the expression for $x$ and using  \eqref{eq:h4Gparam}, we deduce
\begin{equation*}
C=\frac{x}{1+x^2}\ , \qquad G= \frac{x \left(1+x+x^2+x^3+x^4\right)^2}{(1+x)^4 \left(1+x^2\right)^3}\ ,
\end{equation*}
which allows us to now use $x$ instead of $C$ to parametrize $G$ (the parametrization being univocal by requiring $0\leq x \leq 1$,
assuming again $G$ in the range $0\leq G\leq 25/128$).
Rewriting the above expressions for $\alpha$ and $\beta$ in terms of $x$ as 
\begin{equation*}
\alpha = - \frac{1+x+x^2}{(1+x^2)^2}\ , \qquad 
\beta  = - x^2\ \frac{1+x+x^2}{(1+x^2)^2}\ , \\
\end{equation*}
and using $Y_1=-(C+1)=
-(1+x+x^2)/(1+x^2)$, we deduce
\begin{equation*}
Y_k=-\frac{\left(1+x+x^2\right)}{\left(1+x^2\right)^2}\  \frac{\left(1-x^{k+3}\right)}{\left(1-x^{k+1}\right)}\ .
\end{equation*}
Plugging this expression in \eqref{eq:tY} gives
\begin{equation*}
t_k=\frac{x\, \left(1+x+x^2\right)}{\left(1+x+x^2+x^3+x^4\right) }\  \frac{ \left(1-x^{k-1}\right) \left(1-x^{k+5}\right)}{\left(1-x^{k+1}\right)
   \left(1-x^{k+3}\right)}
\end{equation*}
and eventually
\begin{equation*}
r_k=t_k+1=\frac{(1+x)^2 \left(1+x^2\right)}{\left(1+x+x^2+x^3+x^4\right)} \frac{\left(1-x^k\right) \left(1-x^{k+4}\right)}{
   \left(1-x^{k+1}\right) \left(1-x^{k+3}\right)}\ ,
\end{equation*}
with in particular
\begin{equation*}
r_\infty\equiv \lim_{k\to \infty} r_k= \frac{(1+x)^2 \left(1+x^2\right)}{\left(1+x+x^2+x^3+x^4\right)}\ .
\end{equation*}
To eventually find the solution of our original recursion relation \eqref{eq:newrecur}, we simply need to know the value of $R_1$. To this end, we 
combine the identities $R_\infty= R_1 r_\infty$, $R_\infty=1+2 g R_\infty^2$ and $G= g\, R_1^2$ to write
\begin{equation*}
R_1=\frac{1+2G\, r_\infty^2}{r_\infty}=\frac{(1+x+x^2+x^3+x^4)}{\left(1+x^2\right)^2}
\end{equation*}
and finally 
\begin{equation}
\begin{split}
T_k&=R_1\, t_k= 
 \frac{x\, \left(1+x+x^2\right)}{\left(1+x^2\right)^2}\  \frac{ \left(1-x^{k-1}\right) \left(1-x^{k+5}\right)}{\left(1-x^{k+1}\right)
   \left(1-x^{k+3}\right)}\ ,\\
R_k& =R_1\, r_k=\frac{(1+x)^2}{ \left(1+x^2\right)} \frac{\left(1-x^k\right) \left(1-x^{k+4}\right)}{
   \left(1-x^{k+1}\right) \left(1-x^{k+3}\right)}\ , \qquad g=\frac{G}{R_1^2}= \frac{x(1+x^2)}{(1+x)^4}\ . \\
   \end{split}
   \label{eq:solfin}
\end{equation}
We recover here the result guessed in \cite{GEOD} and later recovered in \cite{AB11}. As for the distance-dependent two-point function, it reads:
\begin{equation*}
\hskip -1.cm G_k=R_k-R_{k-1}-\delta_{k,1}=\frac{(1-x)^3 (1+x)^2 \left(1+x+x^2\right)x^{k-1}  \left(1-x^{2 k+3}\right)}{\left(1+x^2\right) \left(1-x^k\right)
   \left(1-x^{k+1}\right) \left(1-x^{k+2}\right) \left(1-x^{k+3}\right)}-\delta_{k,1}
\end{equation*}
for $k\geq 1$.

\section{Hull perimeter statistics}
\label{sec:hull}
Having computed the distance-dependent two-point function $G_k$, we may now evaluate a more refined generating function $H_k(\alpha,d)$ which, in the 
maps enumerated by $G_k$, also controls the \emph{hull perimeter at distance $d$}, namely the length of a particular closed curve separating
the origin $v_0$ of the map and the marked edge of type $(k-1,k)$ and sitting at distance $d$ from $v_0$ (see below for a precise definition of this curve). More precisely, 
if we denote denote by $\mathcal{L}(d)$ the hull perimeter at distance $d$, the generating function $H_k(\alpha,d)$
will now incorporate an extra weight $\alpha^{\mathcal{L}(d)}$ for each map (we have in particular $H_k(1,d)=G_k$ for all the allowed values of $d$). Let us now
come to precise definitions.  
\subsection{Definition of the hull perimeter}
\label{sec:perimeter}
Consider a pointed planar Eulerian triangulation with a marked oriented edge of type $(k-1,k)$, as enumerated by $G_k$, for some $k\geq 3$
and the corresponding $k$-slice, whose apex corresponds to the origin of the Eulerian triangulation, and whose base corresponds
to the marked oriented edge. The dividing line, as we defined it in Section \ref{sec:decomp}, constitutes a particular path 
separating the apex of the slice (which lies strictly above the line) from its base (whose endpoint lies strictly below the line). 
This paths sits moreover ``at distance $d=k-1$" from the apex, in the sense that it is made of a sequence of alternating backward-oriented and forward-oriented edges of type $(d-1,d)$
with $d=k-1$.
In particular, all the vertices lying strictly below the line (i.e.\ on the same side as the base) are at distance at least 
$d$ from the apex. 

Consider now the  \emph{separating line} which simply consists of the dividing line itself minus its first edge $x_0\to y_0$.  If we re-glue the boundary of the slice so as to recover the associated pointed planar Eulerian triangulation, the extremities of this separating line, namely the vertices $y_0$ and $y_p$, 
coalesce so that the line eventually forms a simple closed curve in the pointed map 
which separates the origin of the map from its marked edge and sits at distance $d=k-1$ from the origin, in the above sense. 
The length $2p$ of the separating line defines the \emph{hull perimeter at distance $d$}, here for $d=k-1$, and will be denoted by 
$\mathcal{L}(d)$ as announced\footnote{The part of the map lying on the same side of the (closed) separating line forms what is called the \emph{hull
at distance $d$}, which is the connected domain obtained by removing from the map a particular connected 
component which (i) has all its vertices at distance larger than or equal to $d$ from the origin, and (ii) contains the marked oriented edge.}.
For $k=2$ ($d=1$), the separating line may be viewed as reduced to a single vertex at the apex of the slice (or, after re-gluing, at the origin of the pointed map)
and we set $\mathcal{L}(1)=0$ accordingly.

\begin{figure}
\begin{center}
\includegraphics[width=6cm]{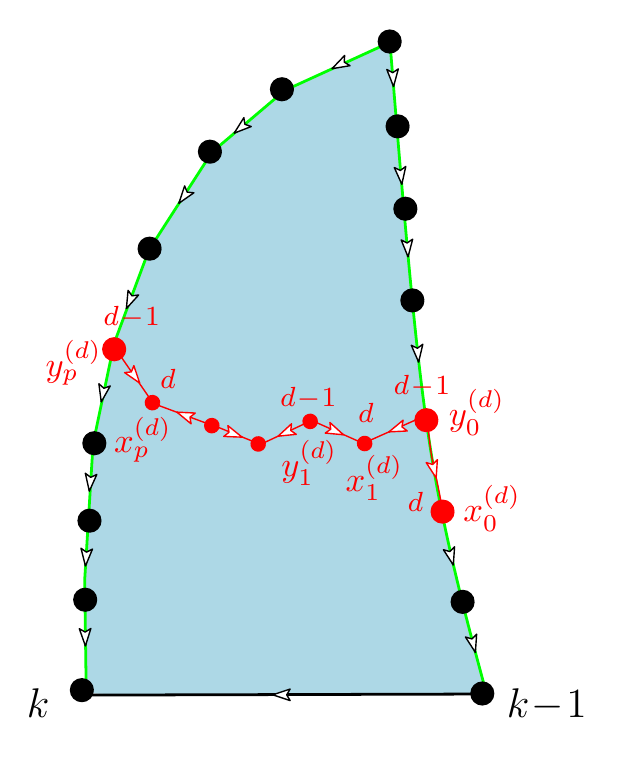}
\end{center}
\caption{Construction of the dividing line at some arbitrary distance $d$
(see text -- the light-blue background is supposedly filled with black and white faces).}
\label{fig:dividingtwo}
\end{figure}
The above definition for $d=k-1$ is easily generalized to some arbitrary $d$ in the range $2\leq d \leq k-1$. As explained in \cite{G16}, a dividing line 
at distance $d$ may be constructed in a $k$-slice as follows (see figure \ref{fig:dividingtwo}): we first pick the unique backward 
oriented edge connecting the right boundary vertices $x^{(d)}_0$ and $y^{(d)}_0$ at respective distance $d$ and $d-1$ 
from the apex and the follow the sequence of leftmost two-step paths $y^{(d)}_m\to x^{(d)}_{m+1} \to y^{(d)}_{m+1}$,
$m\geq 0$ made of vertices satisfying $d(x^{(d)}_{m+1})=d$, $d(y^{(d)}_{m+1})=d-1$, and $x^{(d)}_{m+1}\neq x^{(d)}_{m}$, $y^{(d)}_{m+1}\neq y^{(d)}_{m}$.
With the same arguments as in Section \ref{sec:decomp}, it is easily
shown that this line indeed exists and forms a simple path connecting the right boundary to the left boundary, reached after a number $p^{(d)}$ of steps.
Upon re-gluing the $k$-slice into a pointed Eulerian triangulation, the separating line, consisting now of this new dividing line minus
its first edge $x^{(d)}_0\to y^{(d)}_0$, forms a simple closed curve whose length $\mathcal{L}(d)=2p^{(d)}$ defines the hull perimeter at distance $d$, now for 
some arbitrary $d$ in the range $2\leq d \leq k-1$. 
Again for $d=1$ and arbitrary $k\geq 2$, we view the separating line 
as being reduced to the origin vertex and set $\mathcal{L}(1)=0$.

\subsection{Generating functions with a control on the hull perimeter} 
\label{sec:genhp}
Having defined the hull perimeter $\mathcal{L}(d)$, we now wish to compute the generating function $H_k(\alpha,d)$ incorporating the weight $\alpha^{\mathcal{L}(d)}$.
For starters, consider $\ell$-slices with $2\leq \ell \leq k$, as enumerated by $T_k$, for $k\geq 3$, and the corresponding hull perimeter 
$\mathcal{L}(\ell-1)$ at $d=\ell-1$. To enumerate these $\ell$-slices with a weight $\alpha^{\mathcal{L}(\ell-1)}$,
we simply have, in the slice decomposition of Section \ref{sec:decomp}, to assign a weight $\alpha^2$ to each of the sub-slices constituting the upper part
of the slice. Indeed, each such sub-slice is canonically associated to $2$ edges of the separating line (see figure \ref{fig:dividingone}). At the level of
generating functions, this corresponds to replacing in our recursion the generating function $T_{k-1}$ for each of these sub-slices by $\alpha^2 T_{k-1}$.
In other words, the enumeration of $\ell$-slices with $2\leq \ell \leq k$ and a weight $\alpha^{\mathcal{L}(\ell-1)}$ is achieved by
the quantity
\begin{equation*}
\sum_{p\geq 0} \mathcal{K}_p\, (\alpha^2 T_{k-1})^p= \mathcal{K}(\alpha^2 T_{k-1})\ .
\end{equation*} 

\begin{figure}
\begin{center}
\includegraphics[width=9cm]{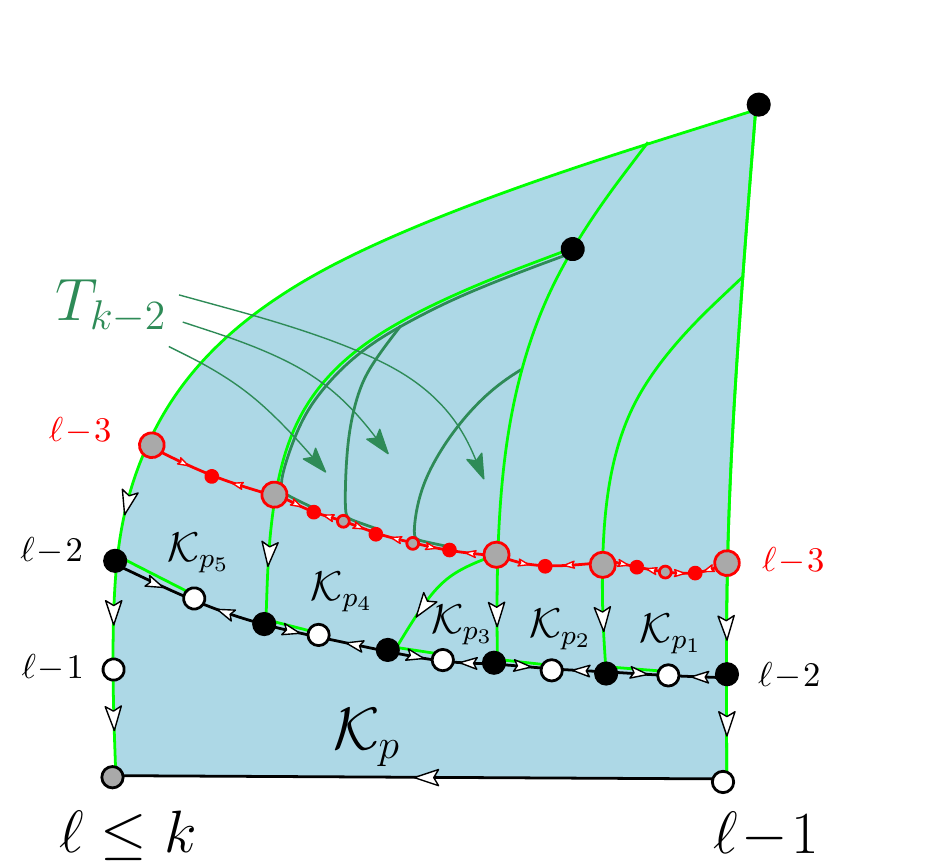}
\end{center}
\caption{Construction of the separating line at distance $d=\ell-2$ (i.e.\ made of vertices at alternating distances $\ell-3$ and $\ell-2$) in an $\ell$-slice 
enumerated by $T_k$ by concatenation of the separating lines at distance $\ell_i-1$ of the $\ell_i$-slices, $i=1,\cdots,p$ (here $p=5$) appearing
as sub-slices in its recursive decomposition (the light-blue background is supposedly filled with black and white faces). The length
of this line (here in red) is $\mathcal{L}(\ell-2)=2\sum_{i=1}^{p}p_i$ ($=14$ here), $2p_i$ being the length of the separating line of the $i$-th sub-slice
(with for instance $p_4=3$ here). We note that $p_i$ is also the number of sub-sub-slices (enumerated by $T_{k-2}$)
in the upper part of the $i$-th sub-slice, the lower part of the sub-slice being itself enumerated by $\mathcal{K}_{p_i}$.}
\label{fig:hull}
\end{figure} 
Consider now the hull perimeter $\mathcal{L}(\ell-2)$ at $d=\ell-2$. As explained in \cite{G16}, the separating line at distance $d=\ell-2$ in an $\ell$-slice 
enumerated by $T_k$ (and satisfying moreover $\ell\geq 3$) is obtained
by \emph{concatenating the separating lines} at distance $\ell_i-1$ of the $\ell_i$-slices, $i=1,\cdots,p$, appearing
as sub-slices in the recursive decomposition of the slice (see figure \ref{fig:hull} for an illustration). To reconstruct the desired weight 
$\alpha^{\mathcal{L}(\ell-2)}$, these sub-slices, originally enumerated by $T_{k-1}=\mathcal{K}(T_{k-2})$, must now be counted by $\mathcal{K}(\alpha^2 T_{k-2})$.
We eventually deduce that the generating function for $\ell$-slices with $2\leq \ell \leq k$ with an extra weight $\alpha^{\mathcal{L}(\ell-2)}$ whenever $\ell\geq 3$
is $\mathcal{K}(\mathcal{K}(\alpha^2 T_{k-2}))$ .

Repeating the argument recursively, we find that, more generally,
\begin{equation*}
\underbrace{\mathcal{K}\big( \mathcal{K}\big(\cdots \big( \mathcal{K}\big(}_{m\ \hbox{\scriptsize times}}\alpha^2\, T_{k-m}\big)\big)\big)\big)
\end{equation*}
is the generating function of $\ell$-slices with $2\leq \ell\leq k$ and with a weight $\alpha^{\mathcal{L}(\ell-m)}$
whenever $\ell\geq m+1$.\footnote{Note that slices with $\ell\leq m$ are enumerated as in $T_k$, with no dependence on $\alpha$.}

As a direct consequence,
the desired generating function $H_k(\alpha,d)\equiv H_k(\alpha,d,g)$ for planar pointed planar Eulerian triangulations 
with a marked edge of type $(k-1,k)$, with a weight $g$ per white face and a weight $\alpha^{\mathcal{L}(d)}$, is given by
\begin{equation}
H_k(\alpha,d)=\underbrace{ \mathcal{K}\big( \mathcal{K}\big(\cdots \big( \mathcal{K}\big(}_{k-d\ \hbox{\scriptsize times}}\alpha^2\, T_{d}\big)\big)\big)\big)- \underbrace{ \mathcal{K}\big( \mathcal{K}\big(\cdots \big( \mathcal{K}\big(}_{k-d\ \hbox{\scriptsize times}}\alpha^2\, T_{d-1}\big)\big)\big)\big)\ .
\label{eq:Hk}
\end{equation}
for any arbitrary $d$ in the range $2 \leq d\leq k-1$. 

The generating function $H_k(\alpha,d)$ can then be computed explicitly from this formal expression as follows: we note that, from the recursion \eqref{eq:relation}, $T_k$ is obtained 
from $T_d$ by $(k-d)$ successive actions of the hernel $\mathcal{K}$ and the explicit form \eqref{eq:solfin} of $T_k$
therefore allows us to write the relation
\begin{equation}
\begin{split}
\hskip -1.2cm
& \underbrace{ \mathcal{K}\bigg( \mathcal{K}\bigg(\cdots \bigg( \mathcal{K}\bigg(}_{k-d\ \hbox{\scriptsize times}}      
 \frac{x\, \left(1\!+\!x\!+\!x^2\right)}{\left(1\!+\!x^2\right)^2}\  \frac{ \left(1\!-\!x^{d-1}\right) \left(1\!-\!x^{d+5}\right)}{\left(1\!-\!x^{d+1}\right)
   \left(1\!-\!x^{d+3}\right) }  \bigg)\bigg)\bigg)\bigg)\\ & 
  \hskip 6cm =  \frac{x\, \left(1\!+\!x\!+\!x^2\right)}{\left(1\!+\!x^2\right)^2}\  \frac{ \left(1\!-\!x^{k-1}\right) \left(1\!-\!x^{k+5}\right)}{\left(1\!-\!x^{k+1}\right)
   \left(1\!-\!x^{k+3}\right)}\ .\\
   \end{split}
   \label{eq:Kiter}
\end{equation}
From this relation, we immediately deduce\footnote{Indeed, writing $x^k=x^{k-d}\, x^d$, we may consider that, for a given $x$, the relation \eqref{eq:Kiter} involves two independent parameters, 
on one hand the integer $k-d$ (which appears both as a power of $x$ and as the number of iterations) and, on the other hand, the variable $x^d$. 
This latter independent parameter may then be 
replaced on both sides of the equation by some arbitrary quantity, for instance $\lambda x^d$ with some arbitrary $\lambda$. In practice, $\lambda$ must be small enough for the equality to remain valid (see \cite{G16} for details).} the more general identity:
\begin{equation*}
\begin{split}
\hskip -1.2cm
& \underbrace{ \mathcal{K}\bigg( \mathcal{K}\bigg(\cdots \bigg( \mathcal{K}\bigg(}_{k-d\ \hbox{\scriptsize times}}      
 \frac{x\, \left(1\!+\!x\!+\!x^2\right)}{\left(1\!+\!x^2\right)^2}\  \frac{ \left(1\!-\!\lambda\, x^{d-1}\right) \left(1\!-\!\lambda\, x^{d+5}\right)}{\left(1\!-\!\lambda\, x^{d+1}\right)
   \left(1\!-\!\lambda\, x^{d+3}\right) }  \bigg)\bigg)\bigg)\bigg)\\ & 
    \hskip 6.cm =
   \frac{x\, \left(1\!+\!x\!+\!x^2\right)}{\left(1\!+\!x^2\right)^2}\  \frac{ \left(1\!-\!\lambda\, x^{k-1}\right) \left(1\!-\!\lambda\, x^{k+5}\right)}{\left(1\!-\!\lambda\, x^{k+1}\right)
   \left(1\!-\!\lambda\, x^{k+3}\right)}\\
   \end{split}
\end{equation*}
for arbitrary (small enough) $\lambda$. 
We have in particular
\begin{equation*}
\underbrace{ \mathcal{K}\big( \mathcal{K}\big(\cdots \big( \mathcal{K}\big(}_{k-d\ \hbox{\scriptsize times}}\alpha^2\, T_{d}\big)\big)\big)\big)=
\frac{x\, \left(1\!+\!x\!+\!x^2\right)}{\left(1\!+\!x^2\right)^2}\  \frac{ \left(1\!-\!\lambda(\alpha,d)\, x^{k-1}\right) \left(1\!-\!\lambda(\alpha,d)\, x^{k+5}\right)}{\left(1\!-\!\lambda(\alpha,d)\, x^{k+1}\right)
   \left(1\!-\!\lambda(\alpha,d)\, x^{k+3}\right) } 
\end{equation*}
if we take for $\lambda(\alpha,d)$ the solution of
\begin{equation*}
\alpha^2 \frac{ \left(1\!-\!x^{d-1}\right) \left(1\!-\!x^{d+5}\right)}{\left(1\!-\!x^{d+1}\right)
   \left(1\!-\!x^{d+3}\right) }=  \frac{ \left(1\!-\!\lambda(\alpha,d)\, x^{d-1}\right) \left(1\!-\!\lambda(\alpha,d)\, x^{d+5}\right)}{\left(1\!-\!\lambda(\alpha,d)\, x^{d+1}\right)
   \left(1\!-\!\lambda(\alpha,d)\, x^{d+3}\right) } 
\end{equation*}
(this equation has two solutions and we must pick that satisfying $\lambda(1,d)=1$, see \cite{G16} for details). We finally deduce
from \eqref{eq:Hk} the explicit expression
\begin{equation*}
\begin{split}
\hskip -1.2cm H_k(\alpha,d) & =\frac{\left(1\!-\!x^2\right)^2 \left(1\!+\!x\!+\!x^2\right)}{1\!+\!x^2}\times \\
&\qquad \times \frac{x^{k-1} (\lambda(\alpha,d\!-\!1)\!-\!x \lambda(\alpha,d)) \left(1\!-\!\lambda(\alpha,d)\lambda(\alpha,d\!-\!1) x^{2 k+3}  \right)}{\left(1\!-\!\lambda(\alpha,d) x^{k+1} \right) \left(1\!-\!\lambda(\alpha,d) x^{k+3} 
   \right) \left(1\!-\!\lambda(\alpha,d\!-\!1) x^k  \right) \left(1\!-\!\lambda(\alpha,d\!-\!1) x^{k+2} \right)}\ , \\
   \end{split}
\end{equation*}
where $x$ parametrizes $g$ via \eqref{eq:solfin} and with $\lambda(\alpha,d)$ defined as above. This expression, valid for $2\leq d\leq k-1$, is sufficient to explore the
statistics of hull perimeters in large planar Eulerian triangulations. For completeness, recall that $\mathcal{L}(1)=0$ so that $H_k(\alpha,1)=G_k$ for all $k\geq 2$.

\subsection{Statistics of hull perimeters in large Eulerian triangulations}
\label{sec:stathp}
We end this Section by giving a number of explicit results on the statistics of hull perimeters. From now on, we will consider maps with \emph{a fixed size},
i.e.\ we will work in the ensemble of pointed 
planar Eulerian triangulations with a marked oriented edge of type $(k-1,k)$ and \emph{with a fixed
number $F$ of white faces}. The number of such maps is $[g^F]G_k$.
We are actually interested in the so-called local limit of large maps, which corresponds to let $F\to \infty$, keeping (at least at a first stage) the
parameter $k$ finite. We denote by $E_k\left(\left\{\cdot\right\}\right)$ the expectation value of some quantity $\left\{\cdot\right\}$
in this limit. We have for instance 
\begin{equation*}
E_k\left(\alpha^{\mathcal{L}(d)}\right)=\lim_{F\to \infty}\frac{[g^F]H_k(\alpha,d)}{[g^F]G_k}\ .
\end{equation*}
The large $F$ behavior of $[g^F]H_k(\alpha,d)$ or $[g^F]G_k$ is easily obtained from the singular behavior of $H_k(\alpha,d)$ and $G_k$
respectively when $g$ approaches the critical value $g^\star=1/8$ (corresponding to $x=1$ in \eqref{eq:solfin}). Setting $g=g^\star (1-\epsilon^4)$ which, from the relation \eqref{eq:solfin}
between $g$ and $x$, amounts to setting
\begin{equation*}
x=\frac{1-\epsilon}{1+\epsilon}\ ,
\end{equation*}
this singular behavior is obtained straightforwardly from the small $\epsilon$ expansion of $H_k(\alpha,d)$ and $G_k$.
Only even powers of $\epsilon$ may appear since $H_k(\alpha,d)$ and $G_k$ have expressions which are invariant under $x\to 1/x$. 
The first two terms of the expansions are of order $\epsilon^0$ and
$\epsilon^4$ (there happens to be no $\epsilon^2$ term) and correspond
to regular terms in $(g^\star-g)$. The most singular behavior therefore comes from the next term, of order $\epsilon^6$, and corresponds to 
a singularity of the form $(g^\star-g)^{3/2}$. We may thus write
\begin{equation*}
H_k(\alpha,d)\vert_{\rm sing.}\sim \mathfrak{H}_k(\alpha,d)(g^\star-g)^{3/2}\ , \qquad G_k\vert_{\rm sing.}\sim \mathfrak{G}_k(g^\star-g)^{3/2}\ ,
\end{equation*}
and the determination of the coefficients $\mathfrak{H}_k(\alpha,d)$
and $\mathfrak{G}_k$ is a rather easy task from the above expressions for $H_k(\alpha,d)$ and $G_k$.
The desired expectation value above is then simply given by
\begin{equation*}
E_k\left(\alpha^{\mathcal{L}(d)}\right)=\frac{\mathfrak{H}_k(\alpha,d)}{\mathfrak{G}_k}\ .
\end{equation*}
A similar calculation was presented in \cite{G16} in the context of quadrangulations or general triangulations and we invite the reader
to look at this reference for explicit intermediate steps. Here we do not reproduce all the details but we only present the final results 
of this straightforward, although slightly involved calculation.
\vskip .3cm
We distinguish two situations: (i) the case of infinite $k$ and (ii) the case of finite $k$.

\subsubsection*{$\triangleright$ Infinite $k$}
Expressions are slightly simpler if, after having sent $F\to \infty$, we also send $k\to \infty$. Note that the limits are taken in this order and, in particular,
$k$ \emph{does not scale with $F$}. For $F\to \infty$, it is expected that, among all the connected domains at distance larger than $d$ from the origin, 
only one has an infinite size (i.e.\ an infinite number of white faces). Letting $k\to \infty$ imposes that the chosen oriented edge of type $(k-1,k)$ 
actually belongs to this infinite connected component. The hull perimeter  $\mathcal{L}(d)$ then corresponds, so to say, to \emph{the length of a
closed curve sitting at distance $d$ from the origin 
and separating this origin from infinity}. For short, we denote by $E_\infty \left(\left\{\cdot\right\}\right)$ the expectation value $E_k \left(\left\{\cdot\right\}\right)$ 
in the limit $k\to\infty$. We find explicitly, for $d\geq 2$:
\begin{equation*}
\hskip -1.cm  E_\infty \left(\alpha^{\mathcal{L}(d)}\right)=\sqrt{\frac{(d+1)(d+3)(9-\alpha^2)+8\alpha^2}{(d+1)(d+3)(1-\alpha^2)+8\alpha^2} }
   -\sqrt{\frac{d(d+2)(9-\alpha^2)+8\alpha^2}{d(d+2)(1-\alpha^2)+8\alpha^2}}
\end{equation*}
which, by differentiation, yields
\begin{equation*}
E_\infty (\mathcal{L}(d))=\frac{3 \left(d^4+6 d^3+10 d^2+3 d-5\right)}{8 (d+1) (d+2)}\ .
\end{equation*}
Expanding the above expression for $E_\infty \left(\alpha^{\mathcal{L}(d)}\right)$ in powers of $\alpha$, we find equivalently that the probability $p_\infty\left(\mathcal{L}(d)=2p\right)$ that $\mathcal{L}(d)$ equals some even integer $2p$ ($p\geq 1$) reads, for $d\geq 2$,
\begin{equation}
\begin{split}
& p_\infty\left(\mathcal{L}(d)=2p\right)=4\times 3^{1-2 p} \left(\frac{(d-1)^p (d+5)^p}{(d+1)^{p} (d+3)^{p}}-\frac{(d-2)^p (d+4)^p}{d^{p} (d+2)^{p}}\right)\, A(p) \\
& A(p)\equiv  \sum_{q=0}^{p-1}2^q\, {p-1\choose q}\, {2q+1\choose q}\ . \\
\end{split}
\label{eq:pinf}
\end{equation}
For large $d$, $\mathcal{L}(d)$ scales as $d^2$ and we define a rescaled hull perimeter $L(d)$ via
\begin{equation*}
L(d)\equiv\frac{\mathcal{L}(d)}{d^2}\ .
\end{equation*}
The latter has a finite limit for large $d$, namely 
\begin{equation}
\lim_{d\to \infty} E_\infty (L(d))=\frac{3c}{2}
\label{eq:Lmoy}
\end{equation}
with a \emph{scaling factor $c$} whose value is
\begin{equation*}
c=\frac{1}{4}\ .
\end{equation*}
\begin{figure}
\begin{center}
\includegraphics[width=11cm]{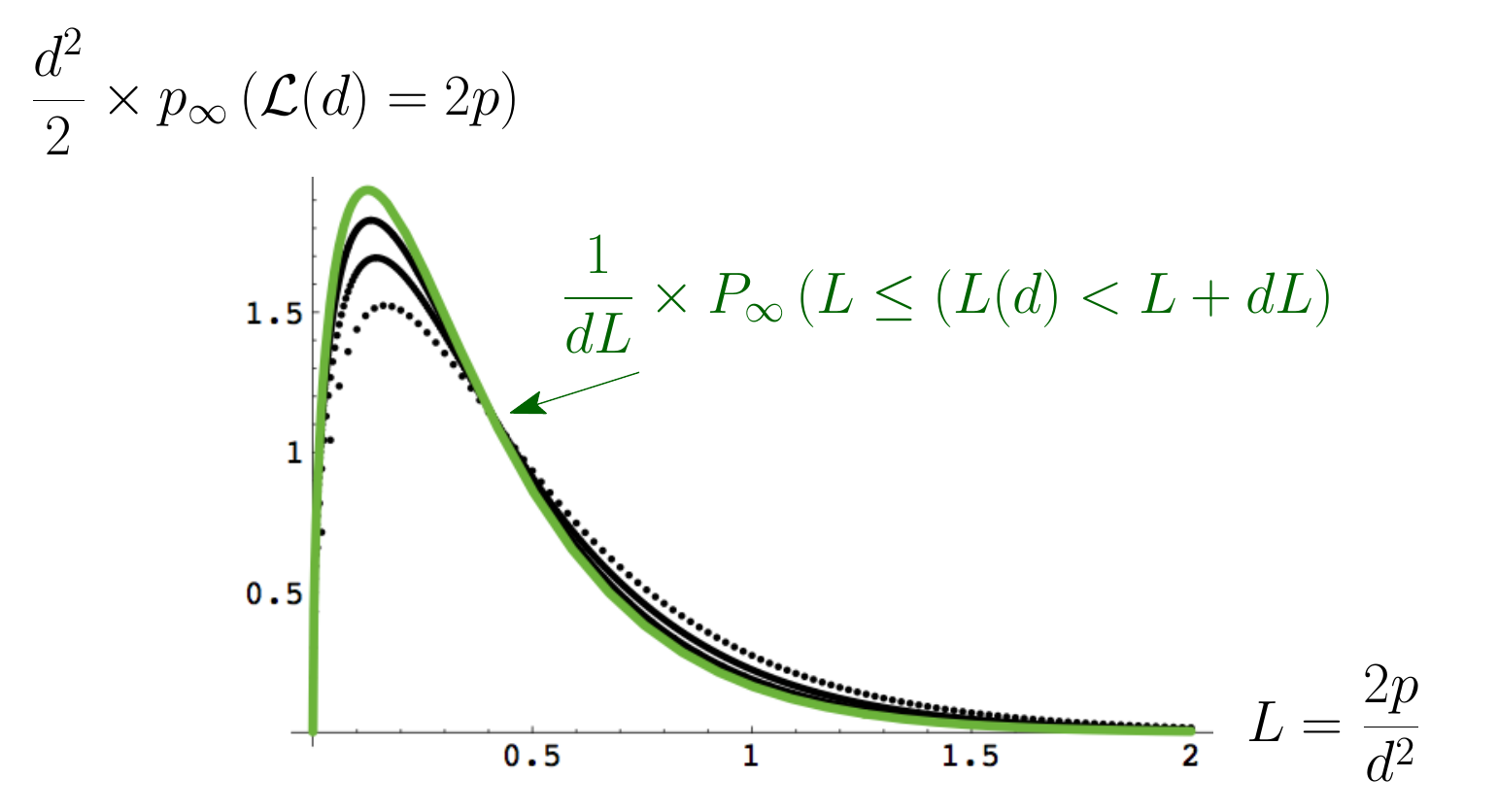}
\end{center}
\caption{The probability $p_\infty\left(\mathcal{L}(d)=2p\right)$ for $d=10$, $20$ and $50$ (with the peak increasing for increasing $d$).
We use as abscissa the quantity $L=2p/d^2$ and we rescale the probability by a factor $d^2/2$ so that it converges to the 
probability density $(1/dL)P_\infty(L\leq L(d)< L+dL)$ with expression \eqref{eq:Pinf}, displayed here as a green solid curve.}
\label{fig:probaperimeter}
\end{figure} 
A more precise characterization of the statistics of $L(d)$ at large $d$ is via the quantity 
\begin{equation}
\lim_{d\to \infty} E_\infty (e^{-\tau L(d)})= \frac{1}{(1+c\, \tau)^{3/2}}
\label{eq:statL}
\end{equation}
or, equivalently, via the probability $P_\infty(L\leq L(d)< L+dL)$ that $L(d)$ lies in the range $L\leq L(d)< L+dL$
(as obtained either by a simple inverse Laplace transform of \eqref{eq:statL} or as the limit of \eqref{eq:pinf}):
\begin{equation}
\lim_{d\to \infty} P_\infty(L\leq L(d)< L+dL)= \frac{2}{\sqrt{\pi}}\frac{\sqrt{L}}{c^{3/2}}e^{-\frac{L}{c}}dL\ .
\label{eq:Pinf}
\end{equation}
The probability $p_\infty\left(\mathcal{L}(d)=2p\right)$ for $d=10$, $20$, $50$ and the corresponding limiting probability density 
when $d\to \infty$ are plotted in figure \ref{fig:probaperimeter}. 
Expressions \eqref{eq:Lmoy}, \eqref{eq:statL} and \eqref{eq:Pinf} are exactly the same as those found in \cite{G16} and in earlier works by Krikun \cite{Krikun03,Krikun05} and by   
Curien and Le Gall \cite{CLG14a,CLG14b}
for other families of maps. Only the value of the scaling factor $c$ (here $c=1/4$) is specific to Eulerian triangulations.
This is a manifestation of the well-known \emph{universality} found in the limit of large maps and for large distances.

\subsubsection*{$\triangleright$ Finite $k$}
Expressions in this case are more involved and we only give here the expectation value of $\mathcal{L}(d)$. We find precisely:
\begin{equation*}
\begin{split}
& \hskip -1.2cm E_k(\mathcal{L}(d))=
\frac{k (k\!+\!1) (k\!+\!2) (k\!+\!3)}{2 (2 k\!+\!3) \left(10 k^6\!+\!90 k^5\!+\!283 k^4\!+\!348 k^3\!+\!103 k^2\!-\!42 k\!-\!36\right)}   \times \\
&  \hskip -1.2cm\times \Bigg(
(d\!-\!1) (d\!+\!1) (d\!+\!3) (d\!+\!5) (k\!+\!2)\times \\ & \times \frac{(k\!+\!1)^2(k\!+\!3)^2 \left(5 k^2\!+\!20 k\!+\!4\right) \!-\!(d\!-\!1) (d\!+\!1) (d\!+\!3) (d\!+\!5) \left(5 d^2\!+\!20 d\!+\!24\right)\!-\!18}{(d\!+\!2) (k\!+\!1)^2 (k\!+\!3)^2}\\
   & \hskip -.8cm-
(d\!-\!2) d (d\!+\!2) (d\!+\!4) (k\!+\!1) \times \\  & \times \frac{k^2 (k\!+\!2)^2 \left(5 k^2\!+\!10 k\!-\!11\right) \!-\!(d\!-\!2) d (d\!+\!2) (d\!+\!4) \left(5 d^2\!+\!10 d\!+\!9\right)\!-\!18}{(d\!+\!1) k^2 (k\!+\!2)^2}   \Bigg)\\\
   \end{split}
\end{equation*}
for $2\leq d\leq k-1$. 
As in \cite{G16}, we may eventually look at the limit where $d$ and $k$ become large simultaneously, keeping the ratio 
$u\equiv d/k$ finite (with in particular $0<u<1$). The rescaled hull perimeter $L(d)$ has a finite expectation value
in this limit, which depends on $u$ only, namely
\begin{equation*} 
\lim_{k\to \infty} E_k (L(k\, u)) = \frac{3c}{2}(1+u-3 u^6
+u^7) \ .\\
\end{equation*}
We recover here, as expected, the universal expression found in \cite{G16}.

\section{Conclusion}
\label{sec:conclusion}

The present study raises a number of natural questions: how general is our recursive approach? Can it be extended to 
other families of maps than triangulations, quadrangulations and Eulerian triangulations?  A common feature of these three cases 
is that their slice generating functions  can be written as bi-ratios of a quantity  $u_k\equiv 1-x^k$, involving a single parameter $x$ 
(parametrizing the weight given to the faces). In our formalism, these bi-ratios themselves come from
the expression of $Y_k$ as a simple ratio of $u_k$'s and this latter structure reveals the existence of some underlying homographic 
recursion relation for $Y_k$ (or equivalently of some underlying geometric recursion for $W_k$, which eventually explains the $x^k$ dependence of $u_k$).
For more general families of maps, the slice generating functions are now bi-ratios of more
involved $u_k$'s whose form is typical of so-called \emph{discrete solitons} (see \cite{GEOD} for details). A first clue
to help generalizing our method would be to understand what type of recursion leads to such a discrete soliton form.

Another natural question concerns the intermediate quantity $Y$ in our calculations. As we noted earlier (see footnote \ref{foot3}), many inequivalent choices of the relation 
between $t$ and $Y$ are in fact possible, which transform our recursion into a simple homographic recursion and eventually lead to the same expression for 
$T_k$.  We may wonder whether some particular choices are more natural than others and, in particular, lead to some $Y_k$ with a direct combinatorial interpretation.

As a natural extension of our result, we note that an expression for the slice generating function of planar Eulerian triangulations was given in \cite{FG14}
which incorporates more parameters by assigning \emph{different weights to the vertices of each (gray, black or white) color}. Again this expression 
is the result of some educated guess and no constructive derivation was provided.
It is easy to incorporate such color-dependent vertex weights in our approach. This then leads to three copies of the generating functions 
$\Phi$ and $\Omega$, determined by three independent closed systems, each depending on the three vertex weights at hand. 
Unfortunately, although there is no fundamental obstacle in using our method to solve these systems, expressions become rather involved and we were not able 
to recover the expression of \cite{FG14}.

Finally, we note that Eulerian triangulations, as $3$-constellations, have an underlying three-fold symmetry corresponding to a cyclic permutation of the colors of their vertices. 
Our construction of the dividing line explicitly breaks this symmetry and involves maps, as those enumerated by $h_{2i}$, where colors 
play inequivalent roles\footnote{This explicit breaking of symmetry is also visible in the rescaling $h_{2i}(g)\to \tilde{h}_{2i}(G)$ which, as mentioned in the footnote \ref{foot2}, 
consists in fact in going from general Eulerian triangulations 
to Eulerian triangulations having no multiple edges connecting a black to a white vertex or a white to a gray vertex, while multiple edges connecting a gray to a black vertex remain allowed.}.
One may wonder if the three-fold symmetry could be preserved in our approach by performing some other decomposition of the slices, involving some other,
more symmetric choice of the dividing line.
We were not able to find such a manageable symmetry-preserving scheme.

\section*{Acknowledgements} 
The author acknowledges the support of the grant ANR-14-CE25-0014 (ANR GRAAL).

\appendix
\section{Derivation of \eqref{eq:recurold}}
\begin{figure}
\begin{center}
\includegraphics[width=13cm]{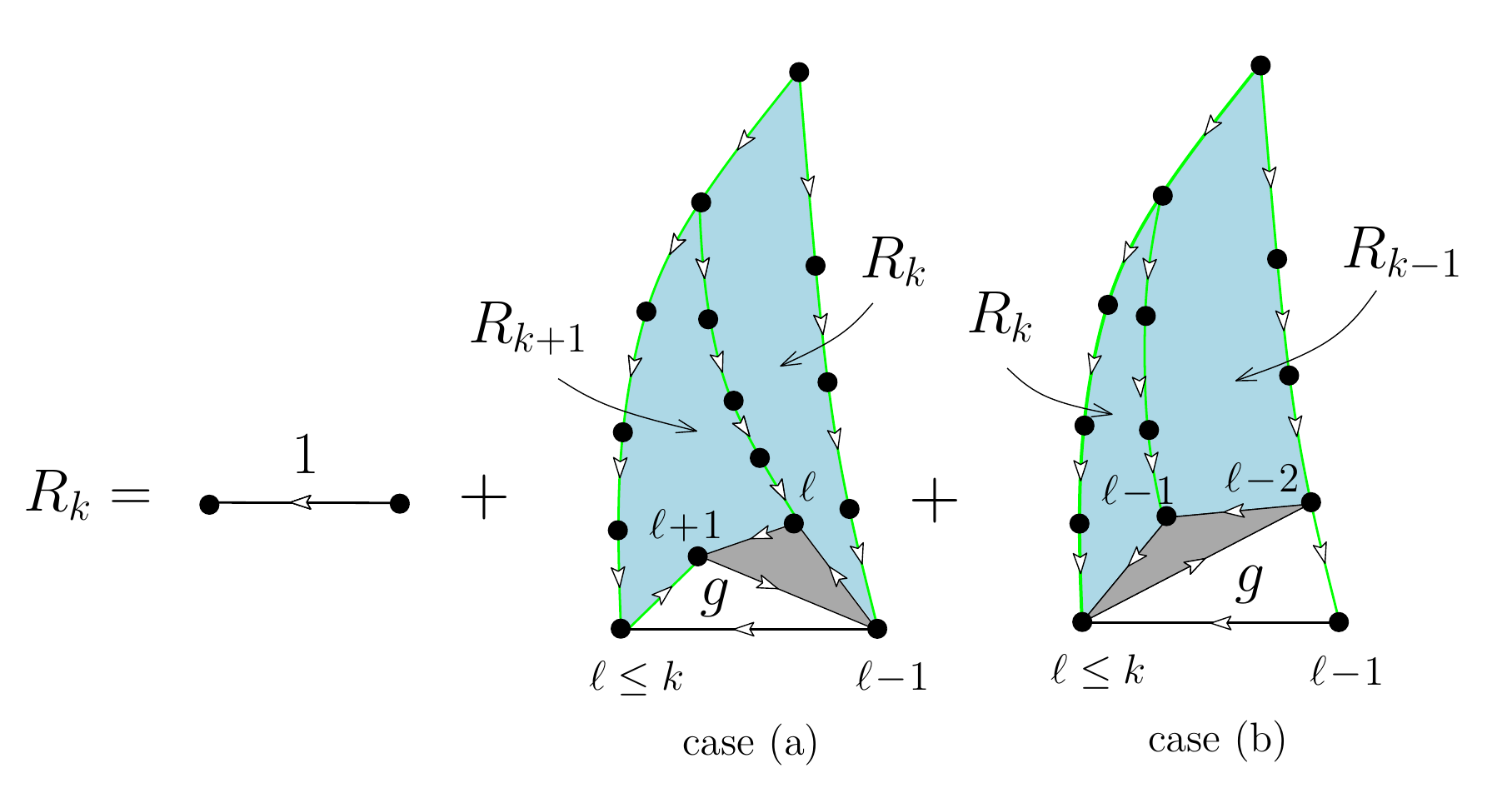}
\end{center}
\caption{A schematic picture explaining the relation \eqref{eq:recurold} for the generating function $R_k$ (the light-blue backgrounds are supposedly filled with black and white faces).}
\label{fig:recur}
\end{figure}
Let us recall how to derive the relation \eqref{eq:recurold}.
By definition, $R_k$ enumerates $\ell$-slices with $1 \leq \ell \leq k$. The first term $1$ in \eqref{eq:recurold} is simply
the contribution of the single-edge-map. For the other $\ell$-slices, we consider 
the white face immediately on the right of the base. By construction, the base in an $\ell$-slice is an oriented edge of type $(\ell-1,\ell)$,
hence a short-edge. As illustrated in figure \ref{fig:recur}, the long-edge incident to the considered white face is either (a) incident to the origin of the 
base (at distance $\ell-1$) hence of type $(\ell+1,\ell-1)$ (the third incident edge being of type $(\ell,\ell+1)$),
or (b), if $\ell\geq 2$, incident to the endpoint of the base (at distance $\ell$) hence of type $(\ell,\ell-2)$ (the third incident edge being of type $(\ell-2,\ell-1)$).
In the first case (a), the white face is of type $\ell-1$ and so is the black face on the other side of the long-edge
(note that long-edges cannot belong to the boundary of the slice which is made of short-edges only, thus both sides
of a long-edge are inner faces). In particular, the third vertex incident to this black face (i.e.\ that which is
not an extremity of the long-edge) is at distance $\ell$ from the apex. Drawing the leftmost
backward shortest path from this vertex to the apex\footnote{The notion of leftmost path usually assumes implicitly that we start from an oriented edge, not from a vertex.
Still it is also well-defined if, starting from a vertex, we forbid a given sector incident to this vertex. Here we forbid the sector 
made of the black face that we just singled out.} and cutting along this line divides the $\ell$-slice into two
parts. These parts are easily seen to be slices, more precisely an arbitrary $(\ell'+1)$-slice, $\ell'\geq 0$ (whose base is the short-edge of
type $(\ell+1,\ell)$ incident to the black face) and an arbitrary $\ell''$-slice, $\ell''\geq 1$ (whose base is the short-edge of
type $(\ell,\ell-1)$ incident to the black face) with $\max(\ell',\ell'')=\ell$.\footnote{\label{foot4}The fact that one of the quantities $\ell'$ 
and $\ell''$ may be less the $\ell$ is due to the fact that the leftmost backward shortest path at which we split the slice hits in general one of the boundaries
of the slice before reaching the apex.} Demanding that $\ell\leq k$ is equivalent to demanding that $1\leq (\ell'+1)\leq k+1$ and 
$1\leq \ell''\leq k$, hence the configurations of case (a) are enumerated by $g\, R_{k+1}R_{k}$, with a weight $g$ for the white face.
This explains the first contribution to the second term in the right hand side of  \eqref{eq:recurold}. In the second case (b), the short-edge of type $(\ell-2,\ell-1)$
incident to the white face is necessarily a right boundary edge as otherwise, this right boundary would not be 
the unique shortest path between the apex and the first extremity of the base. The black face incident to the long-edge
is now of type $\ell-2$ and its third incident vertex is at distance $\ell-1$ from the apex.  Drawing the leftmost
backward shortest path from this vertex to the apex and cutting along this line divides the $\ell$-slice into an arbitrary $\ell'$-slice, $\ell'\geq 1$ (whose base 
is the short-edge of type $(\ell,\ell-1)$ incident to the black face) and an arbitrary $(\ell''-1)$-slice, $\ell''\geq 2$ (whose base is the short-edge of
type $(\ell-1,\ell-2)$ incident to the black face) with $\max(\ell',\ell'')=\ell$.\textsuperscript{\ref{foot4}} Demanding that $\ell\leq k$ is equivalent to demanding that $1\leq \ell'\leq k$ and 
$1\leq \ell''-1\leq k-1$. The configurations in case (b) are thus enumerated by $g\, R_{k}R_{k-1}$ (with the convention $R_0=0$), hence
the second contribution to the second term in the right hand side of \eqref{eq:recurold}.

\bibliographystyle{plain}
\bibliography{eulerianhull}

\end{document}